\documentclass[final]{siamart}

\usepackage{braket,amsfonts}
\usepackage{amsmath, amssymb, amsbsy}
\usepackage{array,tabu,multirow}
\usepackage{graphicx,epstopdf}
\usepackage[caption=false]{subfig}
\usepackage{blkarray}

\newcommand{\pd}[2]{\frac{\partial#1}{\partial#2}}
\renewcommand{\div}{\nabla\cdot}
\newcommand{\grad}{\nabla\,}
\newcommand{\rodt}{\boldsymbol{D}}
\newcommand{\ve}[1]{\boldsymbol{#1}}
\newcommand{\cv}[1]{\boldsymbol{#1}}
\newcommand{\dv}[1]{\mathbf{#1}}
\newcommand{\veu}{\boldsymbol{u}}
\newcommand{\vev}{\boldsymbol{v}}
\newcommand{\vef}{\boldsymbol{f}}
\newcommand{\veg}{\boldsymbol{g}}

\newcommand{\TheTitle}{Preconditioners for Two-Phase Incompressible Navier--Stokes Flow} 
\newcommand{\TheAuthors}{N. Bootland, A. Bentley, C. Kees, and A. Wathen}

\headers{Preconditioners for Two-Phase Flow}{\TheAuthors}

\title{{\TheTitle}\thanks{
This publication is based on work supported by the EPSRC Centre for Doctoral Training in Industrially Focused Mathematical Modelling (EP/L015803/1) in collaboration with the US Army Coastal and Hydraulics Laboratory and HR Wallingford. Permission was granted by the Chief of Engineers to publish this information.}}%

\author{Niall Bootland%
  \thanks{Mathematical Institute, University of Oxford, Radcliffe Observatory Quarter, Woodstock Road, Oxford, OX2 6GG, UK (\email{bootland@maths.ox.ac.uk}, \email{wathen@maths.ox.ac.uk}).}%
  \and
  Alistair Bentley%
  \thanks{Department of Mathematical Sciences, Clemson University, Martin Hall, Clemson, SC 29634, USA (\email{ abentle@clemson.edu}).}%
  \and
  Christopher Kees%
  \thanks{Coastal and Hydraulics Laboratory, US Army Engineer Research and Development Center, 3909 Halls Ferry Road, Vicksburg, MS 39180-6133, USA (\email{christopher.e.kees@usace.army.mil}).}%
  \and
  \mbox{Andrew Wathen}\footnotemark[2]%
}

\Crefname{ALC@unique}{Line}{Lines}

\ifpdf
\hypersetup{
  pdftitle={\TheTitle},
  pdfauthor={\TheAuthors}
}
\fi

\begin{document}

\maketitle

\begin{abstract}
  We consider iterative methods for solving the linearised Navier--Stokes equations arising from two-phase flow problems and the efficient preconditioning of such systems when using mixed finite element methods. Our target application is simulation within the Proteus toolkit; in particular, we will give results for a dynamic dam-break problem in 2D. We focus on a preconditioner motivated by approximate commutators which has proved effective, displaying mesh-independent convergence for the constant coefficient single-phase Navier--Stokes equations. This approach is known as the ``pressure convection--diffusion'' (PCD) preconditioner [H. C. Elman, D. J. Silvester and A. J. Wathen, \emph{Finite Elements and Fast Iterative Solvers: with Applications in Incompressible Fluid Dynamics}, second ed., Oxford University Press, 2014]. However, the original technique fails to give comparable performance in its given form when applied to variable coefficient Navier--Stokes systems such as those arising in two-phase flow models. Here we develop a generalisation of this preconditioner appropriate for two-phase flow, requiring a new form for PCD. We omit considerations of boundary conditions to focus on the key features of two-phase flow. Before considering our target application, we present numerical results within the controlled setting of a simplified problem using a variety of different mixed elements. We compare these results with those for a straightforward extension to another commutator-based method known as the ``least-squares commutator'' (LSC) preconditioner, a technique also discussed in the aforementioned reference. We demonstrate that favourable properties of the original PCD and LSC preconditioners (without boundary adjustments) are retained with the new preconditioners in the two-phase situation.
\end{abstract}



\section{Introduction}
\label{sec:intro}

The motivation for this work stems from the computational challenges involved in simulating two-phase flows, in particular air--water flows. The primary costs incurred in such simulations are in the solution of the linear systems arising after linearisation and discretisation of the governing equations. As such, the efficient solution of these linear systems is crucial and motivates the development of preconditioners appropriate for the variable coefficient nature of two-phase flow.

We consider incompressible flow of two immiscible Newtonian phases. Suppose the two-phase problem is defined on an open bounded domain $\Omega \subset \mathbb{R}^{d}$ (\mbox{$d=2,3$}), with one phase occupying $\Omega_{1}$ and the second $\Omega_{2}$ such that \mbox{$\overline{\Omega} = \overline{\Omega}_{1} \cup \overline{\Omega}_{2}$} and \mbox{$\Omega_{1} \cap \Omega_{2} = \emptyset$}. Let the interface between the two phases be denoted \mbox{$\Gamma = \partial \Omega_{1} \cap \partial \Omega_{2}$}. Note that $\Omega_{1}$ and $\Omega_{2}$ may vary over time and that neither set is required to be connected.

The fluid flow is modelled by the incompressible Navier--Stokes equations
\begin{align}
\label{Navier--Stokes1}
\alpha \rho \, \pd{\veu}{t} + \rho \, \veu \cdot \grad \veu + \grad p - \div \left( 2 \mu \, \rodt \veu \right) & = \rho \vef & \text{in} \ \Omega_{i}, & \\
\label{Navier--Stokes2}
\div \veu & = 0 & \text{in} \ \Omega_{i}, & \ i = 1,2.
\end{align}
Here we have fluid velocity $\veu$, pressure $p$, density $\rho$, dynamic viscosity $\mu$, and body force per unit mass $\vef$. We use the notation \mbox{$\rodt \veu = \frac{1}{2} \left( \grad \veu + (\grad \veu)^{T} \right)$} for the rate of deformation tensor and \mbox{$\cv{\sigma}(\veu,p) = -p \, \cv{I} + 2\mu \, \rodt\veu$} for the stress tensor. The value $\alpha = 0$ is taken for steady flow whilst $\alpha = 1$ for time-dependent flow, in which case the dependent variables may vary over time.

After suitable scaling we assume the variables are dimensionless with the piecewise constant density and viscosity being given by
{\setlength\extrarowheight{9pt}
\begin{align}
\label{DensityandViscosityFunctions}
\rho = \left\{
\begin{array}{ll}
1 & \text{in } \Omega_{1}\\[9pt]
\dfrac{\rho_{2}}{\rho_{1}} & \text{in } \Omega_{2}\\[9pt]
\end{array}
\right., & & \mu = \left\{
\begin{array}{ll}
\dfrac{1}{Re} & \text{in } \Omega_{1}\\[9pt]
\dfrac{\mu_{2}}{\mu_{1}Re} & \text{in } \Omega_{2}\\[9pt]
\end{array}
\right.,
\end{align}
where $\rho_{i}$ and $\mu_{i}$ are, respectively, the (dimensional) density and viscosity of the fluid in $\Omega_{i}$, $i = 1,2$ and $Re$ is the Reynolds number of the first fluid. We are interested in the behaviour for increasing Reynolds numbers and so assume the Reynolds number of the second fluid is smaller than that of the first so that $\rho_{2} / \rho_{1} \le \mu_{2} / \mu_{1}$. The difficulty in solving the problem typically increases for larger Reynolds numbers and so here $Re$ characterises the Reynolds number dependence. 
}

In the time-dependent case we require a technique to keep track of the interface as it moves in time, typically this might be a volume-of-fluid (VOF) method or a level set method (see e.g.\ \cite{Rudman,GrossReicheltAndReusken,KeesEtAl,AkkermanEtAl}). Often, for numerical solution, the tracking of the fluids can be decoupled from solving the Navier--Stokes system so that, when we are required to solve this system, the interface position is known. Since our focus is on the efficient solution of the Navier--Stokes system we suppose that the interface position is known and thus the density and viscosity functions are fully specified. Note that, more generally, these need not be piecewise constant functions, for instance when using a model which incorporates a narrow transition region around the interface.

To complete our model we require appropriate conditions on the boundary $\partial \Omega$ and interface $\Gamma$. Here, for simplicity, we consider enclosed flow in 2D. We suppose that the boundary is split into two parts $\partial \Omega_{D}$ and $\partial \Omega_{FS}$ with $\partial \Omega = \partial \Omega_{D} \cup \partial \Omega_{FS}$ and $\partial \Omega_{D} \cap \partial \Omega_{FS} = \emptyset$. We impose the conditions
\begin{align}
\label{BoundaryConditionsD}
\veu & = \veu_{D} & \text{on} \ & \partial \Omega_{D}, \\
\label{BoundaryConditionsN}
\veu \cdot \cv{n} = 0, \ \left(\cv{\sigma}(\veu,p) \cv{n}\right)\cdot\cv{t} & = 0 & \text{on} \ & \partial \Omega_{FS},\\
\label{BoundaryConditionsGamma}
\left[ \veu \right] = \cv{0}, \ \left[ \cv{\sigma}(\veu,p) \cv{n} \right] & = \veg & \text{on} \ & \Gamma,
\end{align}
where $\cv{n}$ and $\cv{t}$ are unit normal and tangent vectors to the corresponding boundary, \mbox{$\left[a \right]_{\Gamma} = \left( a|_{\Omega_{1}} - a|_{\Omega_{2}} \right)|_{\Gamma}$}, and $\veg$ is a localised force term, for instance corresponding to interface tension in the continuum surface force (CSF) model \cite{BrackbillKotheAndZemach,GrossReicheltAndReusken}. The boundary conditions are interpreted as fully specifying the velocity $\veu_{D}$ on the Dirichlet boundary $\partial \Omega_{D}$ and applying free-slip on the wall boundary $\partial \Omega_{FS}$. Note that for enclosed flow an additional constraint is required to fix the level of the pressure such as specifying the average pressure (see \eqref{FiniteElementSpaces}). Across the interface we enforce that the velocity is continuous and if $\veg = \cv{0}$, so that we neglect surface tension (this will be true in our case), then the second condition ensures that the normal stress is continuous. Finally, in the time-dependent case we also suppose that we have suitable initial conditions.

To numerically solve the problem we apply a nonlinear iteration method to treat the nonlinearity $\rho \, \veu \cdot \grad \veu$ in the momentum equation \eqref{Navier--Stokes1}; in our exposition we consider Picard iteration, though Newton iteration is possible too. In the time-dependent case we suppose that implicit time integration is applied, for instance using the backward Euler scheme. See, for example, \cite{SilvesterElmanKayAndWathen} for further details in the case of single-phase Navier--Stokes flow. The linearised equations that follow are known as the generalised Oseen problem: given a divergence-free vector field $\cv{w}$, find $\veu$ and $p$ such that
\begin{align}
\label{Oseen1}
\frac{\alpha \rho}{\Delta t} \veu + \rho \, \cv{w} \cdot \grad \veu + \grad p - \div \left( 2 \mu \, \rodt \veu \right) & = \cv{r} & \text{in} \ \Omega_{i}, & \\
\label{Oseen2}
\div \veu & = 0 & \text{in} \ \Omega_{i}, & \ i = 1,2,
\end{align}
along with the boundary conditions \eqref{BoundaryConditionsD}--\eqref{BoundaryConditionsGamma}. In Picard iteration the wind $\cv{w}$ is the approximation of $\veu$ from the previous iteration, or else zero in the starting iteration, in which case the corresponding Stokes problem is to be solved. Here $\cv{r}$ collects the forcing term $\rho \cv{f}$ and all other known terms such as those arising from the time integration. Solution of the linearised system \eqref{Oseen1}--\eqref{Oseen2} along with appropriate boundary conditions is key to the overall efficiency of the computational modelling and is the focus for our development of preconditioners.

To discretise in space we consider finite element methodology, nonetheless the preconditioning techniques we will describe could equally be applied to alternative discretisations such as the MAC scheme \cite{HarlowAndWelch}. The details of the weak formulation and finite element discretisation are known and are found in, for example, \cite{GrossReicheltAndReusken}; we summarise the resulting finite element problem here. We assume appropriate, though not necessarily inf--sup stable, finite element spaces for velocity and pressure given by
\begin{align}
\label{FiniteElementSpaces}
\mathbb{V}^{h} \subset H^{1}(\Omega)^{2} & & \text{and} & & \mathbb{Q}^{h} \subset \left\lbrace p \in L^{2}(\Omega) \ : \ \int_{\Omega} p = 0 \right\rbrace
\end{align}
respectively. Then we wish to find $\veu_{h} \in \mathbb{V}^{h}$ and $p_{h} \in \mathbb{Q}^{h}$ such that
\begin{align}
\label{FiniteElementProblem1}
\frac{\alpha}{\Delta t} m(\veu_{h},\vev_{h}) + n(\cv{w}_{h};\veu_{h},\vev_{h}) + b(p_{h},\vev_{h}) + a(\veu_{h},\vev_{h}) & = (\cv{r},\vev_{h}), \\[1ex]
\label{FiniteElementProblem2}
b(q_{h},\veu_{h}) & = c(p_{h},q_{h}),
\end{align}
for all $\vev_{h} \in \mathbb{V}^{h}_{0}$ and $q_{h} \in \mathbb{Q}^{h}$. Here $\cv{w}_{h}$ is the known previous iterate for $\veu_{h}$ as part of the nonlinear iteration and $\mathbb{V}^{h}_{0}$ is the finite element space with homogeneous essential (Dirichlet) boundary conditions applied. Further, $(\cdot,\cdot)$ denotes the standard $L^{2}$ inner product while
\begin{align}
\label{formmandn}
m(\veu,\vev) & = \int_{\Omega} \rho \, \veu \cdot \vev, & n(\cv{z};\veu,\vev) & = \int_{\Omega} \rho \left( \cv{z} \cdot \grad \veu \right) \cdot \vev, \\
\label{formsaandb}
a(\veu,\vev) & = \int_{\Omega} 2 \mu \, \mathrm{Tr} \left( \rodt \veu \, \rodt \vev \right), & b(p,\veu) & = - \int_{\Omega} p \, \div \veu,
\end{align}
and $c(\cdot,\cdot)$ is a stabilisation term, needed for inf--sup unstable spaces, or else is zero. Note that the trilinear form $n(\cdot;\cdot,\cdot)$ may additionally require some stabilisation in the case of dominating advection, for instance by using a streamline-diffusion method (see \cite{SilvesterElmanKayAndWathen} for discussion in the single-phase case).

To obtain a linear system we let $\left\lbrace \cv{\varphi}_{j} \right\rbrace$ be a set of velocity basis functions and $\left\lbrace \psi_{j} \right\rbrace$ be a set of pressure basis functions. Defining the matrices
\begin{align}
\label{Matrices}
\begin{gathered}
M^{(\rho)}_{i,j} = m(\cv{\varphi}_{j},\cv{\varphi}_{i}), \quad N^{(\rho)}_{i,j} = n(\cv{w}_{h};\cv{\varphi}_{j},\cv{\varphi}_{i}), \\[1ex] A^{(\mu)}_{i,j} = a(\cv{\varphi}_{j},\cv{\varphi}_{i}), \quad B_{i,j} = b(\psi_{i},\cv{\varphi}_{j}), \quad C_{i,j} = c(\psi_{j},\psi_{i}),
\end{gathered}
\end{align}
then the generic form of the linear systems to be solved is
\begin{align}
\label{SaddlePointSystem}
\left(
\begin{array}{cc}
\frac{\alpha}{\Delta t} M^{(\rho)} + N^{(\rho)} + A^{(\mu)} & B^{T} \\
B & -C
\end{array}
\right) \left(
\begin{array}{c}
\dv{u} \\
\dv{p}
\end{array}
\right) = \left(
\begin{array}{c}
\dv{f} \\
\dv{g}
\end{array}
\right).
\end{align}
Here the right-hand side vectors $\dv{f}$ and $\dv{g}$ collect together the terms involving $\cv{r}$ and known boundary terms. For clarity, we will use the notation of bracketed exponents when we wish to explicitly illustrate the scaling by density and viscosity within the matrices. The stabilisation matrix $C$ will depend on the choice of finite elements used, as such we defer precise details until discussing such choices in \Cref{sec:MixedFiniteElementParticulars}; note that $C$ might also incorporate some scaling, typically with the inverse of viscosity.

The block system \eqref{SaddlePointSystem} is of (generalised) saddle-point form for the coefficient vectors $\dv{u}$ and $\dv{p}$. For large problems, this system must be solved using iterative methods, typically a Krylov subspace method. Since the presence of $N^{(\rho)}$ makes the system nonsymmetric, a common choice of the Krylov method is GMRES \cite{SaadAndSchultz}; this is what we shall use for our numerical results. The performance of an iterative method depends on the conditioning of the linear system and thus preconditioners are essential for efficient computation. The preconditioners we build on and develop here are block preconditioners which show mesh-independent convergence and only mild dependence on the Reynolds number. Before describing block preconditioning techniques, we discuss related work in preconditioning of the Navier--Stokes equations and similar problems in two-phase flow.

In addition to block preconditioners, a wide variety of preconditioners have been proposed for saddle-point systems originating from incompressible flow and related problems. These include domain decomposition methods (see \cite{QuarteroniAndValli}), however standard approaches may not give robust and scalable results as detailed in \cite{CyrShadidAndTuminaro}. Nonetheless, by incorporating a multilevel approach along with aggressive coarsening based on graph partitioning, as explored in \cite{LinEtAl}, scalability is seen in the tests of \cite{CyrShadidAndTuminaro}. Multigrid techniques are often used for sub-problems in the solution process, however specialist multigrid methods have been proposed for the full Navier--Stokes system going back to Vanka, \cite{Vanka}. While the approach of Vanka uses a coupled smoother, more recently uncoupled smoothers have been seen to offer advantages in efficiency; for example, a method for the Stokes equations using a distributive Gauss--Seidel relaxation based on the least-squares commutator is introduced in \cite{WangAndChen}. Another approach is to use incomplete LU factorisation (ILU), a technique developed in the saddle-point ILU (SILU) preconditioner of \cite{RehmanVuikAndSegal} and the ILU preconditioner for nonsymmetric saddle-point matrices of \cite{KonshinOlshanskiiAndVassilevski}. Augmented Lagrangian techniques can also yield effective methods for the Navier--Stokes equations, in particular see the modified augmented Lagrangian preconditioners of \cite{BenziOlshanskiiAndWang}.

In this work our focus is on block preconditioners, which have seen considerable attention in recent years. Here, a matrix factorisation of the block system is utilised and appropriate approximations of the factors are needed to devise preconditioners which are efficient. The principal challenge is an effective approximation to the Schur complement arising in the factorisation. A taxonomy of these approximate block factorisation (ABF) methods for incompressible Navier--Stokes flow is given in \cite{ElmanEtAl-taxonomy}. Though not originally envisaged within this framework, this includes SIMPLE-type methods \cite{RehmanVuikAndSegal-SIMPLE}. Other popular block preconditioners for Navier--Stokes flow are based on approximate commutators, primarily these are the pressure convection--diffusion (PCD) and least-squares commutator (LSC) preconditioners, which are discussed in detail in \cite{ElmanSilvesterAndWathen}. Though most often described through an approximate commutation relationship, the PCD approach was initially devised in \cite{KayLoghinAndWathen} by considering Green's tensors. The PCD preconditioner requires a convection--diffusion operator projected onto the discrete pressure space and this additional operator must be constructed. In search of a more automatic approach, the LSC preconditioner described in \cite{ElmanHowleShadidShuttleworthAndTuminaro}, but developed from \cite{Elman}, is based solely on algebraic considerations of minimising the norm of the commutator in a least-squares sense. Both of these commutator-based preconditioners show scalable results in the tests of \cite{ElmanSilvesterAndWathen} and \cite{CyrShadidAndTuminaro}.

The preconditioners we have described above are primarily considered in the case of constant density and viscosity flows and, to the authors' knowledge, little has been explored for preconditioning variable density and viscosity Navier--Stokes problems. Nonetheless, recently approaches based on augmented Lagrangian techniques \cite{BenziOlshanskiiAndWang} were developed for variable viscosity in \cite{HeAndNeytcheva} and extended to incompressible non-Newtonian flows in \cite{HeNeytchevaAndVuik}; see also \cite{AxelssonHeAndNeytcheva} for a simpler approach incorporating variable density in the time-dependent case. In the case of variable coefficient problems, such as two-phase flow, most work is devoted to the Stokes problem. Here a pressure mass matrix scaled by the inverse of the viscosity is a good choice of Schur complement preconditioner for the stationary problem, as shown in the case of two-phase flow in \cite{OlshanskiiAndReusken} and more generally investigated for variable viscosity in \cite{GrinevichAndOlshanskii}; see also \cite{BursteddeEtAl,MayBrownAndLePourhiet}. By considering an abstract parameter dependent saddle-point system, this is extended to two-phase non-stationary Stokes flow in \cite{OlshanskiiPetersAndReusken} and can be seen as generalising the Cahouet--Chabard preconditioner \cite{CahouetAndChabard}. Similar block preconditioners are constructed for the variable density and viscosity Stokes problem when using finite volumes in \cite{CaiEtAl}, and applied to uniform staggered grids, showing promising results. This work brings together finite element literature on solving the coupled velocity--pressure system making use of Schur complements and applies such techniques within the finite volume setting, where the dominant paradigm has been splitting or projection methods \cite{GuermondMinevAndShen}. A variant of LSC which takes into account viscosity contrast is given in \cite{MayAndMoresi} and used for their studies in computational geodynamics. Another approach is the Schur method described in \cite{RehmanGeenenVuikSegalAndMacLachlan}; this paper also compares what is effectively the same LSC variant of \cite{MayAndMoresi} and labels this LSC$_{\text{D}}$, a name we shall adopt here too. While these methods hold potential to be used for variable coefficient Navier--Stokes flow we have found no reference to such application. We consider the adaptation of such block preconditioners to two-phase Navier--Stokes flow, in particular the PCD as well as LSC approaches.

The outline for the remainder of this document is as follows. In \Cref{sec:blockpreconditioning} we detail block preconditioning of the saddle-point system before outlining the PCD approach for single-phase flow. Our main contribution is given in \Cref{sec:two-phasepreconditioning} where we extend this preconditioner to two-phase flow. We also describe the LSC approach and provide a simple extension that improves performance in the case of two-phase flow. Numerical results for a simplified problem are given in \Cref{sec:numericalresults} for the choice of \mbox{$\ve{Q}_{2}$--$\,\ve{Q}_{1}$}, \mbox{$\ve{Q}_{1}$--$\,\ve{Q}_{1}$}, and \mbox{$\ve{Q}_{2}$--$\,\ve{P}_{-1}$} elements. For this we use the software package IFISS \cite{IFISS,ElmanRamageAndSilvester} within MATLAB, which we have adapted to incorporate varying density and viscosity. This allows us to describe dependence on the parameters in a controlled manner; we will observe that the favourable properties of the preconditioners (without boundary adjustments) are retained. In \Cref{sec:numericalresults2} we focus on our target application of a dam-break problem simulated within Proteus, a
computational methods
and simulation toolkit
(\texttt{http://proteustoolkit.org}). Since Proteus uses unstructured meshes and \mbox{$\ve{P}_{1}$--$\,\ve{P}_{1}$} elements, we first provide results for our simplified problem in Proteus, verifying similar behaviour to our computations in IFISS, before moving to consider the performance of PCD for the dam-break simulation. Finally, \Cref{sec:conclusions} draws together our conclusions.

\section{Block preconditioning techniques}
\label{sec:blockpreconditioning}

The saddle-point system \eqref{SaddlePointSystem} has a coefficient matrix of the form
\begin{align}
\label{SaddlePointMatrixSystem}
K = \left(
\begin{array}{cc}
F & B^{T} \\
B & -C
\end{array}
\right),
\end{align}
where $F = A^{(\mu)} + N^{(\rho)} + \frac{\alpha}{\Delta t} M^{(\rho)}$. A standard approach of \cite{MurphyGolubAndWathen,Ipsen} for preconditioning, based on a block $LU$-decomposition of $K$, is to right precondition using a matrix of the block triangular form
\begin{align}
\label{BlockPrecon}
P = \left(
\begin{array}{cc}
\widehat{F} & B^{T} \\
0 & - \widehat{S}
\end{array}
\right),
\end{align}
where $\widehat{F}$ and $\widehat{S}$ approximate $F$ and the Schur complement $S=BF^{-1}B^{T}+C$, arising in the decomposition, respectively. If the approximations were exact then only two iterations would be required by the Krylov method \cite{MurphyGolubAndWathen}. However, $F$ and, in particular, $S$ are impractical to work with. Since in application of the preconditioner \eqref{BlockPrecon} we require the actions of $\widehat{F} \, \! ^{-1}$ and $\widehat{S}^{-1}$, we seek approximations where these actions are efficient yet remaining faithful to the operators they approximate. Typically for incompressible fluid flow problems, the difficulty in using the preconditioner \eqref{BlockPrecon} is in finding a good approximation to $S$; see, for instance, \cite{ElmanSilvesterAndWathen} for further details.

As noted in \cite{KayLoghinAndWathen}, for non-self-adjoint problems which yield nonsymmetric linear systems such as \eqref{SaddlePointMatrixSystem}, the development of preconditioners typically relies on heuristic arguments since eigenvalues alone are not generally descriptive of the convergence behaviour for nonsymmetric Krylov methods (though see \cite{KlawonnAndStarke,LoghinAndWathen}). We will focus on the approach of the pressure convection--diffusion (PCD) preconditioner \cite{KayLoghinAndWathen,SilvesterElmanKayAndWathen} for approximating the Schur complement. This is comprehensively described in \cite{ElmanSilvesterAndWathen}; here we briefly outline the key idea and provide the most basic form.
In the remainder of this section we suppose that the density and viscosity are everywhere constant (and so $\rho \equiv 1$) so that we consider single-phase flow. However, to present these ideas in a transparent way for extension to variable coefficient problems we retain the density and viscosity scaling. Note that, due to our nondimensionalisation, $\mu$ here is identified with the kinematic viscosity.

\subsection{The pressure convection--diffusion preconditioner}
\label{sec:PCD}

The approach of the PCD preconditioner is most often motivated by considering the commutator $\mathcal{E}$ between the divergence operator and the convection--diffusion operator\footnote{For a constant viscosity $\mu$ the viscous term, $- \div \left( 2 \mu \rodt \right)$, reduces to the vector Laplacian using the incompressibility constraint; this is the form most often written for exposition of the PCD approach.}
\begin{align}
\label{ConvDiffOperator}
\mathcal{F} & = - \mu \, \div \grad + \rho \, \cv{w}_{h} \cdot \grad.
\end{align}
We suppose that we can have an analogous operator to $\mathcal{F}$ on the pressure space, denoted $\mathcal{F}_{p}$, so that
\begin{align}
\label{Commutator}
\mathcal{E} & = \div \mathcal{F} - \mathcal{F}_{p} \div \, .
\end{align}
Though $\mathcal{F}_{p}$ is in general not rigorously defined, we suppose the commutator is small in some sense. It is noted in \cite{ElmanSilvesterAndWathen} that $\mathcal{E}$ would be zero if $\cv{w}_{h}$ were a constant and the operators were defined on an unbounded domain. The idea is then to use a discrete version of $\mathcal{E}$ and within it to isolate the Schur complement $S$. Upon equating the discrete commutator to be zero we obtain an approximation to $S$ in terms of finite element matrices. When time-stepping, we can define $\mathcal{F}$ to also include the term $\frac{\alpha\rho}{\Delta t}$ from \eqref{Oseen1}. Since this additional term is a scaling of the identity operator it cancels in the commutator. Thus the same approach can be applied in the time-dependent case.

For the discrete commutator, to correctly scale the discrete operators we require the finite element mass matrices on the pressure and velocity spaces given by
\begin{align}
\label{MassMatrices}
M_{p;i,j} & = \int_{\Omega} \psi_{j} \cdot \psi_{i}, & M_{i,j} & = \int_{\Omega} \cv{\varphi}_{j} \cdot \cv{\varphi}_{i}.
\end{align}
The discrete commutator is then
\begin{align}
\label{DiscreteCommutator}
\mathcal{E}_{h} = \left( M_{p}^{-1} B \right) \left( M^{-1} F \right) - \left( M_{p}^{-1} F_{p}^{} \right) \left( M_{p}^{-1} B \right),
\end{align}
where $F_{p}$ represents the discrete form of $\mathcal{F}_{p}$. On equating $\mathcal{E}_{h}$ to be zero this can be rearranged, after post-multiplication by $F^{-1}B^{T}$, to give the approximation
\begin{align}
\label{SchurCompApprox}
B F^{-1} B^{T} & \approx M_{p}^{} F_{p}^{-1} B M^{-1} B^{T}.
\end{align}

For simplicity we will now assume that the pressure approximation is continuous, though extensions can be made to discontinuous pressure approximation \cite{ElmanSilvesterAndWathen}. We will provide results in \Cref{sec:numericalresults} for $\ve{Q}_{2}^{}$--$\,\ve{P}_{-1}^{}$ elements, a pairing with discontinuous pressure. In the continuous case we can define $F_{p}$ as
\begin{align}
\label{Fp}
F_{p;i,j} & = \int_{\Omega} \mu \, \grad \psi_{j} \cdot \grad \psi_{i} + \int_{\Omega} \rho \left( \cv{w}_{h} \cdot \grad \psi_{j} \right) \cdot \psi_{i} + \frac{\alpha}{\Delta t} \int_{\Omega} \rho \, \psi_{j} \cdot \psi_{i}.
\end{align}
These terms can be written as $F_{p} = A_{p}^{(\mu)} + N_{p}^{(\rho)} + \frac{\alpha}{\Delta t} M_{p}^{(\rho)}$ respectively. Finally, the scaled Laplacian term $B M^{-1} B^{T}$ is replaced by the sparse pressure Laplacian $A_{p}$ where
\begin{align}
\label{PressureAp}
A_{p;i,j} & = \int_{\Omega} \grad \psi_{j} \cdot \grad \psi_{i}.
\end{align}
This yields the PCD approximation to the inverse Schur complement
\begin{align}
\label{PCD}
\widehat{S}^{-1} & = A_{p}^{-1} F_{p}^{} M_{p}^{-1},
\end{align}
using the matrices defined in \eqref{PressureAp}, \eqref{Fp}, and \eqref{MassMatrices} respectively. The choice \eqref{PCD} for $\widehat{S}^{-1}$ is also applicable when a stabilisation matrix $C \neq 0$ is necessary \mbox{\cite[Section 9.2.1]{ElmanSilvesterAndWathen}}. Key features of \eqref{PCD} are that we need only multiply by $F_{p}$ and do not have to invert such a matrix, nevertheless the non-normality of the problem is included within the preconditioner through $F_{p}$. For a practical implementation, multigrid methods can be used to effect the action of $A_{p}^{-1}$ and $M_{p}$ can be replaced by the spectrally equivalent matrix $\mathrm{diag}(M_{p})$ and thus easily inverted \cite{Wathen}. Alternatively, the action of $M_{p}^{-1}$ on a vector is well approximated using Chebyshev semi-iteration with only a small number of iterations required \cite{WathenAndRees}.

We remark that the commutator can also be taken with the gradient instead of the divergence, as originally done in \cite{SilvesterElmanKayAndWathen}. The effect of this is that the ordering of the operators in \eqref{PCD}, our approximate inverse Schur complement $\widehat{S}^{-1}$, is reversed. It was seen when considering boundary conditions that the choice of the divergence was favourable; see \cite[Remark 9.3]{ElmanSilvesterAndWathen}.

\subsection{The relation of PCD with the Cahouet--Chabard preconditioner}
\label{sec:Cahouet-ChabardRelation}

We conclude this section by noting a relation between the PCD preconditioner \eqref{PCD} and the Cahouet--Chabard preconditioner \cite{CahouetAndChabard} for the generalised Stokes problem; this will provide a useful viewpoint on how a two-phase PCD preconditioner should behave.

In Stokes flow the convective term is omitted so that in \eqref{SaddlePointMatrixSystem} for single-phase flow \mbox{$F = \mu A + \frac{1}{\Delta t} M$}, where $A$ is the discrete velocity Laplacian and $M$ is the velocity mass matrix (since $\rho \equiv 1$). Cahouet and Chabard show that the appropriate Schur complement which balances these two terms is
\begin{align}
\label{CahouetChabard}
\widehat{S}^{-1} = \mu M_{p}^{-1} + \frac{1}{\Delta t} A_{p}^{-1},
\end{align}
where $A_{p}$ is a discrete pressure Laplacian operator, for instance given by \eqref{PressureAp} for continuous pressure approximation.

Now consider the PCD preconditioner in this case. Since the convective term is not present we have $F_{p} = \mu A_{p} + \frac{1}{\Delta t} M_{p}$ and thus the Schur complement approximation
\begin{align}
\label{PCD=CC}
\widehat{S}^{-1} = A_{p}^{-1} F_{p}^{} M_{p}^{-1} = A_{p}^{-1} \left( \mu A_{p}^{} + \frac{1}{\Delta t}M_{p}^{} \right) M_{p}^{-1} = \mu M_{p}^{-1} + \frac{1}{\Delta t} A_{p}^{-1}.
\end{align}
Hence we see that the PCD preconditioner precisely reduces to the Cahouet--Chabard preconditioner in the case of generalised Stokes flow. This provides the viewpoint that PCD extends the Cahouet--Chabard preconditioner to the case of Navier--Stokes flow.

\section{Preconditioners for two-phase flow}
\label{sec:two-phasepreconditioning}

We now move to the focus of our work, namely preconditioning for two-phase flow. The primary new feature in the problem is the variable density and viscosity which gives additional scaling in the equations. The key to adapting the preconditioners detailed in \Cref{sec:blockpreconditioning} will be to incorporate this scaling appropriately into their formulation.

\subsection{Two-phase pressure convection--diffusion preconditioning}
\label{sec:Two-phasePCD}

On its first introduction, the PCD preconditioner for steady flow was derived using Fourier techniques and Green's tensors \cite{KayLoghinAndWathen}. Due to the variable coefficient nature of two-phase flow, this approach cannot apply here. However, in the original approach, in order that the preconditioner defaults to the optimal choice in the Stokes limit when the convective term tends to zero, a mass matrix is included to give the correct scaling. The same philosophy applies here.

Before addressing the appropriate choice in the Stokes case, we must be clear on the form of the viscous term used. In the literature different forms are considered, including $-\div(\mu\,\grad\veu)$ in \cite{OlshanskiiAndReusken}, $-\div(\mu\,\rodt\veu)$ in \cite{GrinevichAndOlshanskii,OlshanskiiPetersAndReusken}, and $-\div(2\mu\,\rodt\veu)$ in \cite{GrossReicheltAndReusken,CaiEtAl}. The overall scaling of the latter term contains an additional factor of two compared with the former choices, as such an appropriate Schur complement approximation for the Stokes problem differs by a factor of a half. The distinction between different forms of the viscous term used (including incorporating bulk viscosity) and the effect on the Schur complement is discussed further in \cite{CaiEtAl}. Since we use the latter choice above, we present results for this case, including a factor of two where appropriate as compared with the cited references.

For two-phase flow, the appropriate pressure mass matrix arising in the Stokes case is given in \cite{OlshanskiiAndReusken} and scales inversely with the viscosity as
\begin{align}
\label{PressureMassMatrixInvVisc}
M_{p;i,j}^{(1/\mu)} & = \int_{\Omega} (2\mu)^{-1} \psi_{j} \cdot \psi_{i}.
\end{align}
This mass matrix is also suitable for general variable viscosity Stokes flows \cite{GrinevichAndOlshanskii}.

Looking at the construction of the PCD preconditioner in \Cref{sec:PCD}, $F_{p}$ is given as in \eqref{Fp}, but now with the density and viscosity being piecewise constant. If we assume that PCD still takes the form in \eqref{PCD} but with additional scaling, then, so that the viscosity scaling within the Schur complement is commensurate, when we use the scaled mass matrix $M_{p}^{(1/\mu)}$ we require the pressure Laplacian-type term
\begin{align}
\label{PressureApVisc}
A_{p;i,j}^{(\mu)} & = \int_{\Omega} \mu \, \grad \psi_{j} \cdot \grad \psi_{i}.
\end{align}
Note this is already constructed as part of $F_{p}$ just as in the single-phase case since for constant viscosity \mbox{$A_{p}^{(\mu)} = \mu A_{p}$}. These choices give a Schur complement approximation
\begin{align}
\label{PCDVisc}
\widehat{S}^{-1} & = \left(A_{p}^{(\mu)}\right)^{-1} F_{p}^{} \left(M_{p}^{(1/\mu)}\right)^{-1}.
\end{align}
However, from our numerical experience, it is apparent that the performance of this preconditioner depends poorly on the density ratio of the two fluids. In particular, \eqref{PCDVisc} does not work so well in the time-dependent case. To understand this we go back to the relationship with the Cahouet--Chabard preconditioner detailed in \Cref{sec:Cahouet-ChabardRelation}.

By considering an abstract parameter dependent saddle-point system, the authors of \cite{OlshanskiiPetersAndReusken} derive a preconditioner for two-phase time-dependent Stokes flow which can be seen as generalising the Cahouet--Chabard preconditioner. This preconditioner uses an approximation to the Schur complement given by
\begin{align}
\label{GeneralisedCahouetChabard}
\widehat{S}^{-1} = \left(M_{p}^{(1/\mu)}\right)^{-1} + \frac{1}{\Delta t} \left(A_{p}^{(1/\rho)}\right)^{-1},
\end{align}
where the pressure Laplacian-type term is now inversely proportional to the density,
\begin{align}
\label{PressureApInvDens}
A_{p;i,j}^{(1/\rho)} & = \int_{\Omega} \rho^{-1} \grad \psi_{j} \cdot \grad \psi_{i}.
\end{align}
The new feature here is the dependence on the density $\rho$. This can be understood heuristically by considering the single-phase case as follows. Firstly, we note that the viscosity $\mu$ in \eqref{CahouetChabard} is in fact identified with the kinematic viscosity. In our case, where $\mu$ is the dynamic viscosity, the inverse Schur complement is multiplied by the density (in addition to there being a factor of two arising from the choice of viscous term) and so the approximation should read
\begin{align}
\label{DensityScaledCahouetChabard}
\widehat{S}^{-1} = 2\mu M_{p}^{-1} + \frac{\rho}{\Delta t} A_{p}^{-1} = \left( (2\mu)^{-1} M_{p}\right)^{-1} + \frac{1}{\Delta t} \left( \rho^{-1} A_{p}\right)^{-1}.
\end{align}
In the two-phase case, the scaling with density and viscosity must be incorporated within the integral definition of the matrices and so we obtain \eqref{GeneralisedCahouetChabard}. An inverse Schur complement approximation very similar to \eqref{GeneralisedCahouetChabard}, appropriate for finite volume Stokes solvers, is given in \cite{CaiEtAl}. Here the MAC discretisation on staggered grids is employed. The primary difference is that for finite volumes a diagonal matrix of viscosities at each pressure degree of freedom can be used for the viscous term, whereas in the finite element setting a mass matrix is appropriate.

It is now clear, from the viewpoint of the relation between PCD and the Cahouet--Chabard preconditioner, that the approximation \eqref{PCDVisc} does not default to the correct choice of the preconditioner \eqref{GeneralisedCahouetChabard} in the time-dependent Stokes case. To give a more robust generalisation of PCD to two-phase flow with appropriate scaling we must split the matrix $F_{p}$ and treat the separate terms accordingly.

The matrix $F_{p}$ in the two-phase case consists of the three parts
\begin{align}
\label{two-PhaseFpAsASumOfMatrices}
F_{p} & = A_{p}^{(\mu)} + N_{p}^{(\rho)} + \frac{\alpha}{\Delta t}M_{p}^{(\rho)}.
\end{align}
To treat the viscous part $A_{p}^{(\mu)}$ of $F_{p}$ we use the scaling from \eqref{PressureMassMatrixInvVisc} and \eqref{PressureApVisc} but to treat the remaining part, depending on the density, we scale using the mass matrix
\begin{align}
\label{PressureMassMatrixDens}
M_{p;i,j}^{(\rho)} & = \int_{\Omega} \rho \, \psi_{j} \cdot \psi_{i}
\end{align}
and use the density scaled Laplacian-type term \eqref{PressureApInvDens}. This yields the approximation
\begin{align}
\label{GeneralisedPCD(rho)}
\widehat{S}^{-1} = \left(M_{p}^{(1/\mu)}\right)^{-1} + \left(A_{p}^{(1/\rho)}\right)^{-1} \left( N_{p}^{(\rho)} + \frac{\alpha}{\Delta t} M_{p}^{(\rho)} \right) \left(M_{p}^{(\rho)}\right)^{-1}.
\end{align}
In practice however, we have found it beneficial to further cancel the scaling with the density $\rho$ in the final two bracketed terms of \eqref{GeneralisedPCD(rho)} to give two-phase PCD as
\begin{align}
\label{GeneralisedPCD}
\widehat{S}^{-1} = \left(M_{p}^{(1/\mu)}\right)^{-1} + \left(A_{p}^{(1/\rho)}\right)^{-1} \left( N_{p}^{(1)} + \frac{\alpha}{\Delta t} M_{p}^{(1)} \right) \left(M_{p}^{(1)}\right)^{-1}.
\end{align}
Here the bracketed exponent is given as 1 to be clear that these are unscaled terms whose definition does not include $\rho$. Namely we have
\begin{align}
\label{Np}
N_{p;i,j}^{(1)} & = \int_{\Omega} \left( \cv{w}_{h} \cdot \grad \psi_{j} \right) \cdot \psi_{i},
\end{align}
while $M_{p}^{(1)}$ is the standard mass matrix. While this gives a different approximation, the overall scaling of the two separate terms in the sum remain the same and, although it is feasible to cancel the density as in \eqref{GeneralisedPCD}, it remains important that $A_{p}^{(1/\rho)}$ keeps the correct scaling. Further, note that $M_{p}^{(1)}$ remains the same at each time-step while $M_{p}^{(\rho)}$ would change due to the moving phases and so would need recomputing. From the arguments presented here, it is not clear that \eqref{GeneralisedPCD} should be preferable to \eqref{GeneralisedPCD(rho)}, however, we are actively investigating this further.

As with single-phase PCD, in the two-phase PCD approximation of \eqref{GeneralisedPCD} only one pressure Laplacian-type solve is needed, this is because the appearance of $A_{p}^{(\mu)}$ cancels out. We note that in the time-dependent Stokes case, when $N_{p}^{(1)}$ is zero, we have
\begin{align}
\label{GeneralisedPCD=GeneralisedCC}
\begin{split}
\widehat{S}^{-1} & = \left(M_{p}^{(1/\mu)}\right)^{-1} + \left(A_{p}^{(1/\rho)}\right)^{-1} \left( \frac{1}{\Delta t} M_{p}^{(1)} \right) \left(M_{p}^{(1)}\right)^{-1} \\[1ex]
& = \left(M_{p}^{(1/\mu)}\right)^{-1} + \frac{1}{\Delta t} \left(A_{p}^{(1/\rho)}\right)^{-1},
\end{split}
\end{align}
and so we return to the generalised Cahouet--Chabard preconditioner \eqref{GeneralisedCahouetChabard}. Further, despite the new form of the preconditioner \eqref{GeneralisedPCD}, it is a generalisation of the original PCD preconditioner \eqref{PCD} since for everywhere constant viscosity $\mu$ and density $\rho \equiv 1$ (and thus using the Laplacian for the viscous term, removing the factor two in \eqref{PressureMassMatrixInvVisc})
\begin{align}
\label{GeneralisedPCDSimplifiesToPCD}
\begin{split}
\widehat{S}^{-1} & = \left(\mu^{-1} M_{p}\right)^{-1} + A_{p}^{-1} \left( N_{p}^{} + \frac{\alpha}{\Delta t} M_{p}^{} \right) M_{p}^{-1} \\[1ex]
& = A_{p}^{-1} \left( \mu A_{p}^{} \right) M_{p}^{-1} + A_{p}^{-1} \left( N_{p}^{} + \frac{\alpha}{\Delta t} M_{p}^{} \right) M_{p}^{-1} \\[1ex]
& = A_{p}^{-1} F_{p}^{} M_{p}^{-1}.
\end{split}
\end{align}

We remark that the generalisation of PCD to two-phase flow in \eqref{GeneralisedPCD} cannot be written in the form of the original PCD in \eqref{PCD} with scaling in $A_{p}$ and $M_{p}$. The scaling here comes from the mass matrices on the velocity and pressure space which depend on the corresponding norms (see \cite{OlshanskiiPetersAndReusken} for the choice of norms in the generalised Stokes case). If simple scaled norms are used in the mass matrix scaling of the terms in the discrete commutator \eqref{DiscreteCommutator}, one for the pressure space and one for the velocity space, then it is not possible to construct a preconditioner which defaults to the appropriate generalised Cahouet--Chabard preconditioner. We are currently exploring whether commutator-based approaches might still be extensible to derive an appropriate two-phase preconditioner. Nonetheless, it appears that such a form as \eqref{PCD} is a simplification which can only be made for an everywhere constant density and viscosity flow.

While we have considered discretisation with finite elements, we envisage that the approach can also be utilised with finite volumes or finite differences. For instance, using the MAC discretisation \cite{HarlowAndWelch} and building on the work in \cite{CaiEtAl}. The primary requirement is an appropriate discretisation of the pressure space convection term.

Finally, we note that the features of PCD described in \Cref{sec:PCD} are retained in that we do not have to invert a pressure space convection--diffusion operator and only one discrete Laplacian-type sub-problem must be solved in applying \eqref{GeneralisedPCD}. In practice, each mass matrix can be replaced by a diagonal approximation or Chebyshev semi-iteration can be used to approximate the inverse action of each mass matrix on a vector. Multigrid methods are again applicable for the solution of the scaled pressure Laplacian-type sub-problem.

Before moving on to discuss our numerical results, we first briefly detail another commutator-based preconditioner often considered alongside PCD: the least-squares commutator (LSC) method.


\subsection{The least-squares commutator preconditioner}
\label{sec:Two-phaseLSC}

While still using the commutator \eqref{Commutator}, the LSC preconditioner uses a purely algebraic approach to defining $F_{p}$ by choosing it so that the discrete commutator \eqref{DiscreteCommutator} is small in a least-squares sense. It turns out that it is more convenient to consider the adjoint of the commutator and choose $F_{p}$ by minimising each individual vector norm of the columns of $F_{p}^{T}$; for further details see \cite{ElmanSilvesterAndWathen}. This yields the weighted least-squares problem
\begin{align}
\label{LSCWeightedLeastSquares}
\min \left\| \left[ M^{-1} F^{T} M^{-1} B^{T} \right]_{j} - M^{-1} B^{T} M_{p}^{-1} \left[ F_{p}^{T} \right]_{j} \right\|_{M},
\end{align}
for each column $j$. Using the normal equations, $F_{p}$ can be given as
\begin{align}
\label{LSCFp}
F_{p} & = \left( B M^{-1} F M^{-1} B^{T} \right) \left( B M^{-1} B^{T} \right)^{-1} M_{p},
\end{align}
and hence, substituting this expression into \eqref{SchurCompApprox}, we arrive at the Schur complement approximation
\begin{align}
\label{LSCSchurCompApprox}
B F^{-1} B^{T} & \approx \left( B M^{-1} B^{T} \right) \left( B M^{-1} F M^{-1} B^{T} \right)^{-1} \left( B M^{-1} B^{T} \right).
\end{align}
In practice we replace $M$ with its diagonal $T = \mathrm{diag}(M)$, and so construct the sparse scaled Laplacian $B T^{-1} B^{T}$. The LSC inverse Schur complement approximation then takes the form
\begin{align}
\label{LSC}
\widehat{S}^{-1} & = \left( B T^{-1} B^{T} \right)^{-1} \left( B T^{-1} F T^{-1} B^{T} \right) \left( B T^{-1} B^{T} \right)^{-1}.
\end{align}
This approximation is not immediately applicable to stabilised elements where $C \neq 0$, however appropriate modifications can be made in such an instance \cite{ElmanHowleShadidSilvesterAndTuminaro}. As with the PCD method, the action of the inverse of the discrete Laplacian term can be effectively applied using multigrid methods, however, note that we now need to solve two discrete Laplacian sub-problems in the application of $\widehat{S}^{-1}$. Multiplication by $F$ captures the non-normality of the problem within the action of $\widehat{S}^{-1}$.

While the LSC approach is purely algebraic, incorporating boundary adjustments is important in order to obtain scalable performance for enclosed flow problems such as those we consider; see \cite{ElmanSilvesterAndWathen}. We do not pursue this topic here as we focus on the key aspects of generalising the original preconditioners to be more effective in the case of two-phase flow. Since we expect poorer scalability from LSC in this situation, for brevity we choose to forgo discussion on the stabilisation terms required when using inf--sup unstable elements, which follow in a straightforward manner from \cite{ElmanHowleShadidSilvesterAndTuminaro}, and we will omit the LSC approach when discussing our target dam-break simulation.



The scaling within LSC is known to be an important factor. In \cite{Elman}, a precursor to LSC is derived which has the same form as \eqref{LSC} except that the scaling matrix $T$ is not present. That so-called $BFBt$ method was seen to have poorer performance than the later LSC method \cite{ElmanHowleShadidShuttleworthAndTuminaro}. In general, the choice of scaling is essential for the efficiency of such a preconditioner. In an effort to use a scaling which incorporates changes in viscosity for the Stokes problem, \cite{MayAndMoresi} use the form
\begin{align}
\label{LSC_D}
\widehat{S}^{-1} & = \left( B D^{-1} B^{T} \right)^{-1} \left( B D^{-1} F D^{-1} B^{T} \right) \left( B D^{-1} B^{T} \right)^{-1},
\end{align}
where $D_{i,i} = \max_{j}|F_{i,j}|$ is a diagonal scaling matrix. In \cite{RehmanGeenenVuikSegalAndMacLachlan} the authors consider $D = \mathrm{diag}(F)$, which in practice is essentially the same scaling, and label this method LSC$_\text{D}$. Our experience shows that, for the Navier--Stokes case, this preconditioner does not scale well with the problem size, with it exhibiting a poor $h$-dependence. To improve performance over both LSC$_\text{D}$ and the original (single-phase) approach we consider a scaled mass matrix coming from a viscosity weighted least-squares problem.

The choice of norms on the velocity and pressure spaces determines the form of the mass matrices used in the discrete commutator \eqref{DiscreteCommutator} and weighted least-squares problem \eqref{LSCWeightedLeastSquares}. However, the resulting LSC approximation is independent of the norm on the pressure space. In view of this, we introduce the viscosity scaled norm on the velocity space and the corresponding velocity mass matrix
\begin{align}
\label{ViscNorm}
\|\veu\|_{\mu}^{2} & = \int_{\Omega} \mu \, \veu \cdot \veu, & M^{(\mu)}_{i,j} & = \int_{\Omega} \mu \, \cv{\varphi}_{j} \cdot \cv{\varphi}_{i},
\end{align}
and keep the standard $L^2$ norm on the pressure space. Then the construction of the discrete commutator and weighed least-squares problem is the same as in \eqref{DiscreteCommutator} and \eqref{LSCWeightedLeastSquares}, except with $M$ replaced by $M^{(\mu)}$; as such we omit the details here. Letting \mbox{$T^{(\mu)} = \mathrm{diag} ( M^{(\mu)} )$}, the resulting two-phase LSC Schur complement approximation is
\begin{align}
\label{GeneralisedLSC}
\widehat{S}^{-1} & = \left( B {T^{(\mu)}}^{-1} B^{T} \right)^{-1} \left( B {T^{(\mu)}}^{-1} F {T^{(\mu)}}^{-1} B^{T} \right) \left( B {T^{(\mu)}}^{-1} B^{T} \right)^{-1}.
\end{align}
Unlike with PCD, it is not necessary to split the convection--diffusion term and scale differently. Our experience is that such a splitting of $F$ does not gain any improvement in the preconditioner. We note that when the viscosity $\mu$ is everywhere constant this can be factored out and cancelled to give the original form \eqref{LSC}. As with the original LSC preconditioner, \eqref{GeneralisedLSC} requires the solution of two scaled discrete Laplacian-type sub-problems which, in practice, can be achieved using multigrid techniques.

Unlike the PCD method, the LSC approach does not simplify to the Cahouet--Chabard preconditioner in the Stokes case and so we cannot use this viewpoint to guide our choice of scaling. We will see that our choice of viscosity scaling in \eqref{GeneralisedLSC} can significantly improve the robustness of the LSC approach to contrasts in the material parameters, however, it remains unclear whether this choice is optimal. We are currently investigating the importance of such scaling and the impact of boundary adjustments when LSC-type methods are applied to two-phase flow problems.

\subsection{A comparison of PCD and LSC}
\label{sec:ComparisonOfPCDandLSC}

Detailed comparisons between the original PCD and LSC preconditioners are given in \cite{OlshanskiiAndVassilevski} and \cite{ElmanSilvesterAndWathen}. Our two-phase strategies compare similarly. To summarise the key differences: PCD requires the construction of additional matrices on the pressure space while LSC requires only matrices which are readily available; this means that the LSC method directly applies when discontinuous pressure approximation is used while the formulation of PCD is further complicated. On the other hand, PCD naturally extends to stabilised elements while LSC does not immediately apply; here, stabilisation terms are required \cite{ElmanHowleShadidSilvesterAndTuminaro}. Further, only one solution of a discrete Laplacian-type sub-problem is needed by PCD while LSC requires two. We also note that we have not considered boundary conditions for $A_{p}$ and $F_{p}$ (or $N_{p}$) within PCD. A more in depth look at the commutator can shed light on this area but is beyond the scope of this work; see \cite{ElmanSilvesterAndWathen}. Such considerations do not explicitly arise for the LSC preconditioner but can nonetheless be important.

In the single-phase case, mesh-independent convergence rates have been observed for both preconditioners. Eigenvalue bounds can be found in \cite{OlshanskiiAndVassilevski} and while for PCD they are $h$-independent, for LSC the known bounds depend on $h$. Norm-equivalence of the PCD preconditioner is shown in \cite{LoghinAndWathen} and, further, rigorous GMRES convergence bounds based on the field of values are given for a mass matrix Schur complement approximation; see also bounds in \cite{BenziAndOlshanskii}. These field of values bounds are $h$-independent. For PCD and LSC, all known bounds depend on the Reynolds number. Numerical experiments in \cite{ElmanSilvesterAndWathen} show iteration counts for both methods that are mildly dependent on $Re$. Our results in the next section will demonstrate that the favourable properties of these preconditioners for the single-phase case are retained in the two-phase case.

\section{Numerical results for a simplified problem}
\label{sec:numericalresults}

In this section we present results exhibiting the behaviour of our new preconditioners for a benchmark 2D test problem, that of lid driven cavity flow. We suppose \mbox{$\Omega = (-1, 1)^{2}$} with \mbox{$\Omega_{2} = (-\frac{1}{2}, \frac{1}{2})^{2}$}, and thus \mbox{$\Omega_{1} = \Omega \setminus \overline{\Omega}_{2}$}, and consider a regularised cavity with the prescribed flow
\begin{align}
\label{RegularisedCavity}
\veu = \left(
\begin{array}{c}
1 - x^{4} \\
0
\end{array}
\right) \quad \text{for} \quad y = 1, \ -1 \le x \le 1,
\end{align}
along the lid while no flow boundary conditions are imposed on the remaining three edges. Note that this configuration of fluids is similar to that in the SINKER problem used in \cite{MayAndMoresi,RehmanGeenenVuikSegalAndMacLachlan}. We suppose that no external body forces are acting and neglect surface tension. In the governing linearised equations of \eqref{Oseen1}--\eqref{Oseen2} we suppose for simplicity that $\cv{r} = \cv{0}$ so that in the time-dependent case we start from an initially non-moving flow. Further, we solve the Stokes problem (with $\cv{w} = \cv{0}$) using the generalised Cahouet--Chabard preconditioner \eqref{GeneralisedCahouetChabard} to provide the initial guess for the Picard iteration.

The computations are implemented within MATLAB using the software package IFISS \cite{IFISS,ElmanRamageAndSilvester}, which we have adapted to incorporate spatially varying density and viscosity, though note there is no capability to update the phases in time. As such we consider solving for a single time-step or else the more challenging case of a steady problem. Results for a fully dynamic two-phase problem are presented in \Cref{sec:numericalresults2}.

Numerical results provided in this section tabulate the number of iterations that GMRES requires for solving the linear system \eqref{SaddlePointSystem} to a given tolerance. The GMRES Krylov method is used (without restarting) with the preconditioner \eqref{BlockPrecon} in which the inverse Schur complement approximations of two-phase PCD \eqref{GeneralisedPCD} and two-phase LSC \eqref{GeneralisedLSC} are used. We present average GMRES iteration counts over the course of the Picard iteration, omitting the initial Stokes solve. Since our focus is on the quality of the Schur complement approximation, we suppose that the approximation of $F$ is exact, that is $\widehat{F} \, \! ^{-1} = F^{-1}$ as given by a direct solver (`backslash' in MATLAB).

\subsection{Convergence and termination}
\label{sec:ConvergenceAndTermination}

To measure the convergence we look at the residual vector for the linear system given in \eqref{SaddlePointSystem} to solve for the $k+1$th Picard iterate. The linear residual is
\begin{align}
\label{LinearResidual}
\dv{r}^{(\ell)} = \left(
\begin{array}{cc}
\frac{\alpha}{\Delta t} M^{(\rho)} + N^{(\rho)}(\cv{w}_{h}^{(k)}) + A^{(\mu)} & B^{T} \\
B & 0
\end{array}
\right) \left(
\begin{array}{c}
\dv{u}^{(\ell)} \\
\dv{p}^{(\ell)}
\end{array}
\right) - \left(
\begin{array}{c}
\dv{f} \\
\dv{g}
\end{array}
\right),
\end{align}
where $\dv{u}^{(\ell)}$ and $\dv{p}^{(\ell)}$ are the approximate solutions to the system after $\ell$ iterations of preconditioned GMRES. The notation of $N^{(\rho)}(\cv{w}_{h}^{(k)})$ is used here to recall that the convective part depends on the solution $\cv{w}_{h}^{(k)}$ of the previous Picard iteration. Denoting the discrete velocity and pressure vectors of the $k$th Picard iterate by $\dv{w}^{(k)}$ and $\dv{q}^{(k)}$ respectively, we also define the nonlinear residual as
\begin{align}
\label{NonlinearResidual}
\dv{s}^{(k)} = \left(
\begin{array}{cc}
\frac{\alpha}{\Delta t} M^{(\rho)} + N^{(\rho)}(\cv{w}_{h}^{(k)}) + A^{(\mu)} & B^{T} \\
B & 0
\end{array}
\right) \left(
\begin{array}{c}
\dv{w}^{(k)} \\
\dv{q}^{(k)}
\end{array}
\right) - \left(
\begin{array}{c}
\dv{f} \\
\dv{g}
\end{array}
\right).
\end{align}
After the $k$th Picard iteration, we terminate the GMRES iteration for solving the linearised system given in \eqref{SaddlePointSystem} to give the $k+1$th Picard iterate once the relative residual norm decreases below a prescribed tolerance $\varepsilon$, that is we terminate once
\begin{align}
\label{ConvergenceCriterion}
\|\dv{r}^{(\ell)}\|_{2} \le \varepsilon \|\dv{s}^{(k)}\|_{2}.
\end{align}
For each inner application of GMRES we use an initial guess of the past solution $\dv{u}^{(0)} = \dv{w}^{(k)}$, $\dv{p}^{(0)} = \dv{q}^{(k)}$ and use a tolerance of $\varepsilon = 10^{-6}$ for termination. This method of tracking the convergence is the default used by IFISS. The Picard iteration is terminated upon a relative reduction of $10^{-5}$ for the nonlinear residual, as in \cite{ElmanSilvesterAndWathen}. Further, at each Picard iteration, the implementation is such that the linear system is solved for a correction to the current iterate (see \cite[Section 8.3]{ElmanSilvesterAndWathen}).

\subsection{Choice of mixed finite elements}
\label{sec:MixedFiniteElementParticulars}

For our computations we consider three choices of mixed finite element pairings. In each case, we will utilise a regular grid of square elements having side length $h$, thus $4 h^{-2}$ elements in total. First, we use the inf--sup stable pairing of $\ve{Q}_{2}$--$\,\ve{Q}_{1}$ elements. In this case no stabilisation is needed ($C=0$) and since the pressure approximation is continuous we use \eqref{PressureApInvDens} and \eqref{Np} in the formulation of two-phase PCD.

Secondly, we employ equal order $\ve{Q}_{1}$--$\,\ve{Q}_{1}$ elements. The pressure space remains the same but now stabilisation is required and we follow the approach developed in \cite{DohrmannAndBochev}. As noted in \cite{ElmanHowleShadidSilvesterAndTuminaro}, the lack of stability of $\ve{Q}_{1}$--$\,\ve{Q}_{1}$ approximation stems from the mismatch between the discrete divergence of the velocity field and the discrete pressure space $\ve{Q}_{1}$. The approach of \cite{DohrmannAndBochev} is to project the pressure into the appropriate space. This requires the $L^{2}$ projection from the pressure space into $\ve{P}_{0}$, which we denote by $\Pi_{0}$. It is defined locally on each element $\square$ through a local averaging
\begin{align}
\label{P0Projection}
\Pi_{0} p_{h} \vert_{\square} = \frac{1}{\lvert \square \rvert} \int_{\square} p_{h} \qquad \text{for all} \ \square \in \mathcal{T}_{h},
\end{align}
where $\mathcal{T}_{h}$ is the set of $\ve{Q}_{1}$ elements used and $\lvert \square \rvert$ denotes the area of element $\square$. For the Navier--Stokes problem, to ensure commensurate scaling, the resulting stabilisation matrix is scaled by the inverse of viscosity, as in \cite{ElmanHowleShadidSilvesterAndTuminaro}. With variable viscosity the appropriate stabilisation term is then given by
\begin{align}
\label{PressureStabilisationOperatorVisc}
c(p_{h},q_{h}) = \sum_{\square\in\mathcal{T}_{h}} \int_{\square} \mu^{-1} (p_{h}-\Pi_{0}p_{h})|_{\square}(q_{h}-\Pi_{0}q_{h})|_{\square}.
\end{align}
from which the associated stabilisation matrix $C$ in \eqref{SaddlePointSystem} can be constructed. This matrix can be assembled from element contributions in the standard manner.

Finally, we consider a finite element pairing which uses discontinuous pressure approximation. For this we choose $\ve{Q}_{2}^{}$--$\,\ve{P}_{-1}^{}$ elements, which are naturally inf--sup stable and conserve mass elementwise. Two-phase PCD requires the construction of both a Laplacian-type matrix $A_{p}^{(1/\rho)}$ and a convection matrix $N_{p}^{(1)}$ on the pressure space but now definitions \eqref{PressureApInvDens} and \eqref{Np} no longer apply. To construct appropriate operators when the pressure is discontinuous across elements, we do so locally on elements and utilise edge contributions. The construction is more involved with the main ideas being detailed in \cite[Section 9.2.1]{ElmanSilvesterAndWathen} and a proof of concept implementation being found within IFISS \cite{IFISS}. As such, we do not provide a full description here but note that the only difference in our case is the inverse density scaling within the term $A_{p}^{(1/\rho)}$. For this we have found it beneficial to take an arithmetic average of the two density values in the edge contributions when density differs across an element edge. Hence, for the pressure degree of freedom $p_{0}$ at the centroid of a square element $\square$, the corresponding approximation to the Laplacian-type term $A_{p}^{(1/\rho)}$ is given by (cf. \cite[equations (9.29)--(9.30)]{ElmanSilvesterAndWathen})
\begin{align}
\label{discPressureLaplacian}
\begin{split}
\div \left(\rho^{-1} \grad p\right)\vert_{\square} \approx \quad & \frac{1}{h} \left( \frac{2}{\rho_{+x}+\rho_{0}} \frac{p_{+x}-p_{0}}{h} - \frac{2}{\rho_{0}+\rho_{-x}} \frac{p_{0}-p_{-x}}{h} \right) \\ + & \frac{1}{h} \left( \frac{2}{\rho_{+y}+\rho_{0}} \frac{p_{+y}-p_{0}}{h} - \frac{2}{\rho_{0}+\rho_{-y}} \frac{p_{0}-p_{-y}}{h} \right),
\end{split}
\end{align}
where the subscripts denote values from the neighbouring elements: for example, $p_{+x}$ refers to the pressure degree of freedom at the centroid of the neighbouring element in the positive $x$ direction while $\rho_{+x}$ refers to the density in this element. Here, we have assumed for simplicity that the interface lies along edges of the elements used.

\subsection{Preconditioner performance results}
\label{sec:PreconditionerPerformanceResults}

We first note that the numerical experiments presented in this section have also been conducted using ideal variants of the preconditioners, that is using direct solvers for any Laplacian-type and mass matrix sub-problems. However, since the resulting iteration counts are very similar, differing by at most a small number of iterations, we omit them here and justify the following choices. When requiring the action of $\widehat{S}^{-1}$, the solution of the scaled Laplacian-type sub-problems are given approximately using one algebraic multigrid (AMG) V-cycle using a Ruge and St\"{u}ben implementation, specifically the \texttt{HSL\_MI20} code \cite{BoyleMihajlovicAndScott} built into IFISS. For the PCD method \eqref{GeneralisedPCD}, the scaled mass matrix solutions are approximated using three steps of Chebyshev semi-iteration \cite{WathenAndRees}. We present the majority of our results using $\ve{Q}_{2}$--$\,\ve{Q}_{1}$ elements to demonstrate the behaviour of the preconditioners over differing parameters. We only consider this discretisation for two-phase LSC as additional stabilisation is required when using $\ve{Q}_{1}$--$\,\ve{Q}_{1}$ elements and, further, we note that qualitatively similar results were found for $\ve{Q}_{2}^{}$--$\,\ve{P}_{-1}^{}$ elements and thus these results are omitted for brevity.


\begin{table}[t]
\centering
\footnotesize
\caption{Preconditioned GMRES iterations using two-phase PCD \eqref{GeneralisedPCD} / two-phase LSC \eqref{GeneralisedLSC} for the \mbox{$\ve{Q}_{2}$--$\,\ve{Q}_{1}$} solution to the steady lid driven cavity problem with $Re = 100$, $h = 1 / 128$, and varying density ratio $\widehat{\rho} = \rho_{2} / \rho_{1}$ and viscosity ratio $\widehat{\mu} = \mu_{2} / \mu_{1}$. Values marked ``$\star$'' are omitted since here the dominating Reynolds number changes, being associated with the second phase.}
\label{TableDensityViscosity}
\tabulinesep=1.15mm
\begin{tabu}{|c|c|c|c|c|c|c|c|}
\hline
\multirow{2}{*}{$\widehat{\mu}$} & \multicolumn{7}{c|}{$\widehat{\rho}$} \\
\cline{2-8}
& $10^{-3}$ & $10^{-2}$ & $10^{-1}$ & $1$ & $10$ & $100$ & $1000$ \\
\hline
$10^{-3}$ & 30 / 59 & $\star$ & $\star$ & $\star$ & $\star$ & $\star$ & $\star$ \\
\hline
$10^{-2}$ & 29 / 58 & 30 / 58 & $\star$ & $\star$ & $\star$ & $\star$ & $\star$ \\
\hline
$10^{-1}$ & 24 / 54 & 24 / 54 & 25 / 54 & $\star$ & $\star$ & $\star$ & $\star$ \\
\hline
$1$ & 19 / 38 & 19 / 38 & 19 / 38 & 20 / 40 & $\star$ & $\star$ & $\star$ \\
\hline
$10$ & 24 / 44 & 24 / 44 & 24 / 44 & 24 / 44 & 27 / 44 & $\star$ & $\star$ \\
\hline
$100$ & 26 / 38 & 26 / 38 & 26 / 38 & 27 / 38 & 27 / 38 & 29 / 37 & $\star$ \\
\hline
$1000$ & \hphantom{0}26 / 36\hphantom{0} & \hphantom{0}26 / 36\hphantom{0} & \hphantom{0}26 / 36\hphantom{0} & \hphantom{0}27 / 36\hphantom{0} & \hphantom{0}27 / 36\hphantom{0} & \hphantom{0}27 / 34\hphantom{0} & \hphantom{0}28 / 33\hphantom{0} \\
\hline
\end{tabu}
\end{table}

\begin{table}[t]
\centering
\footnotesize
\caption{Preconditioned GMRES iterations using original (single-phase) PCD \eqref{PCD} / LSC \eqref{LSC} for the \mbox{$\ve{Q}_{2}$--$\,\ve{Q}_{1}$} solution to the steady lid driven cavity problem with $Re = 100$, $h = 1 / 128$, and varying density ratio $\widehat{\rho} = \rho_{2} / \rho_{1}$ and viscosity ratio $\widehat{\mu} = \mu_{2} / \mu_{1}$. Values marked ``$\star$'' are omitted since here the dominating Reynolds number changes, being associated with the second phase.}
\label{TableDensityViscosityOriginals}
\tabulinesep=1.15mm
\begin{tabu}{|c|c|c|c|c|c|c|c|}
\hline
\multirow{2}{*}{$\widehat{\mu}$} & \multicolumn{7}{c|}{$\widehat{\rho}$} \\
\cline{2-8}
& $10^{-3}$ & $10^{-2}$ & $10^{-1}$ & $1$ & $10$ & $100$ & $1000$ \\
\hline
$10^{-3}$ & 173 / 261 & $\star$ & $\star$ & $\star$ & $\star$ & $\star$ & $\star$ \\
\hline
$10^{-2}$ & 104 / 74\hphantom{0} & 101 / 74\hphantom{0} & $\star$ & $\star$ & $\star$ & $\star$ & $\star$ \\
\hline
$10^{-1}$ & 46 / 38 & 46 / 38 & 45 / 38 & $\star$ & $\star$ & $\star$ & $\star$ \\
\hline
$1$ & 29 / 37 & 29 / 37 & 29 / 37 & 30 / 39 & $\star$ & $\star$ & $\star$ \\
\hline
$10$ & 36 / 36 & 36 / 36 & 36 / 36 & 36 / 36 & 43 / 36 & $\star$ & $\star$ \\
\hline
$100$ & 76 / 81 & 76 / 81 & 76 / 81 & 76 / 81 & 76 / 80 & 88 / 71 & $\star$ \\
\hline
$1000$ & 125 / 320 & 125 / 320 & 125 / 320 & 125 / 320 & 123 / 322 & 120 / 310 & 140 / 296 \\
\hline
\end{tabu}
\end{table}

\begin{table}[t]
\centering
\footnotesize
\caption{The total number of degrees of freedom (DOFs) in the lid driven cavity problem for the differing grid sizes $h$ and choice of mixed finite elements. For $\ve{Q}_{2}$--$\,\ve{Q}_{1}$ elements this is given by the formula $2(4h^{-1}-1)^{2}+(2h^{-1}+1)^2$. For $\ve{Q}_{2}$--$\,\ve{P}_{-1}$ elements it is given by $2(4h^{-1}-1)^{2}+3(2h^{-1})^2$ while for $\ve{Q}_{1}$--$\,\ve{Q}_{1}$ elements it is $2(2h^{-1}-1)^{2}+(2h^{-1}+1)^2$. Note that in each row of the table the number of velocity degrees of freedom is the same.}
\label{TableNDOFs}
\tabulinesep=1.15mm
\begin{tabu}{|c|c|c||c|c|}
\hline
$h$ & $\ve{Q}_{2}$--$\,\ve{Q}_{1}$ & $\ve{Q}_{2}$--$\,\ve{P}_{-1}$ & $h$ & $\ve{Q}_{1}$--$\,\ve{Q}_{1}$ \\
\hline
$\frac{1}{16}$ & $9,\!027$ & $11,\!010$ & $\frac{1}{32}$ & $12,\!163$ \\
\hline
$\frac{1}{32}$ & $36,\!483$ & $44,\!546$ & $\frac{1}{64}$ & $48,\!899$ \\
\hline
$\frac{1}{64}$ & $146,\!691$ & $179,\!202$ & $\frac{1}{128}$ & $196,\!099$ \\
\hline
$\frac{1}{128}$ & $588,\!291$ & $718,\!850$ & $\frac{1}{256}$ & $785,\!411$ \\
\hline
$\frac{1}{256}$ & $2,\!356,\!227$ & $2,\!879,\!490$ & $\frac{1}{512}$ & $3,\!143,\!683$ \\
\hline
\end{tabu}
\end{table}

To start, we consider steady flow ($\alpha = 0$) and show how the performance of each preconditioner depends on density and viscosity. \Cref{TableDensityViscosity} displays results for two-phase PCD / two-phase LSC for a variety of density ratios \mbox{$\widehat{\rho} = \rho_{2} / \rho_{1}$} and viscosity ratios \mbox{$\widehat{\mu} = \mu_{2} / \mu_{1}$} with given Reynolds number $Re = 100$ and grid size $h = 1/128$. 
Note that we omit values marked ``$\star$'' since here the dominating Reynolds number changes --- dependence on $Re$ is given in \Cref{TableGridSizeReynoldsNumber}. We see that, while there is some variation with the density and, in particular, viscosity ratios, the performance is fairly robust across the range of ratios tested. This is not true of the original PCD or LSC preconditioners for this problem; these methods give performance which depends significantly on the density and viscosity ratios, and poor mesh-dependence away from ratios close to unity. Equivalent values to those in \Cref{TableDensityViscosity} reach into several hundreds of iterations for the original methods, as exhibited in \Cref{TableDensityViscosityOriginals}. We remark that the discrepancy between the two-phase PCD and original PCD method for the case of $\widehat{\rho} = 1$ and $\widehat{\mu}=1$ is due to the additional factor of two present in the two-phase PCD method, stemming from the choice of viscous term used, as discussed in \Cref{sec:Two-phasePCD}.

We note here that some variation is also seen in the results for the generalised Cahouet--Chabard preconditioner in \cite{OlshanskiiPetersAndReusken}; in particular, for large ratios in opposite directions for density and viscosity the method becomes much less efficient. Results provided in \cite{HeAndNeytcheva}, for an augmented Lagrangian approach to solve the Navier--Stokes equations with variable viscosity, also depend on the viscosity through the minimum and maximum values taken in the domain.

Before continuing to give results for varying grid sizes we detail the scope of the overall problem sizes that we use for the comparisons in this work. For each choice of mixed finite elements we consider a sequence of five grids. The total number of degrees of freedom (unknowns in the problem) for each linear system solve are given in \Cref{TableNDOFs}, including a one-dimensional null space related to a constant pressure mode, thus our constraint on the pressure in \eqref{FiniteElementSpaces} to specify it uniquely. We note that on the finest grids the number of degrees of freedom reaches into several million unknowns.

The dependence of the iteration counts on Reynolds number and grid size is given in \Cref{TableGridSizeReynoldsNumber} for the example of $\widehat{\rho} = 1.2 \times 10^{-3}$ and $\widehat{\mu} = 1.8 \times 10^{-2}$, corresponding to ratios in an air--water system. As for the original PCD and LSC preconditioners, there is a mild dependence on $Re$, with higher Reynolds number flows requiring more iterations. These results are comparable with \cite[Table 9.5]{ElmanSilvesterAndWathen}, where the grid sizes used in \Cref{TableGridSizeReynoldsNumber} correspond to grid parameters $l=6,\ldots,10$ of \cite{ElmanSilvesterAndWathen}. (Note that in \cite{ElmanSilvesterAndWathen} ``old'' refers to the original preconditioners without adjustments for boundary conditions, a topic not pursued here as we choose to focus on the key aspects of generalising PCD and LSC to two-phase flow. Further, the operators in PCD act in reverse order; see \cite[Remark 9.3]{ElmanSilvesterAndWathen}.) Our results for two-phase PCD display nearly mesh-independent behaviour akin to the original method applied to a single fluid, with only a very small increase in iterations required for finer grids. For the highest Reynolds number flows, we note a slight jump in iterations required for the finest grid and this is attributed to the fact that new features of these flows are resolved using this grid; this will also be seen in further results. Two-phase LSC shows more dependence on $h$ but the increase per refinement is not too great and again performance is similar to that of the original method (without boundary adjustments) applied to single-phase flow. We remark that the trends shown in \Cref{TableGridSizeReynoldsNumber} are also seen for other choices of density and viscosity ratios. We further note that when a larger number of iterations is required, such as with high Reynolds number flows, the variation seen in iteration counts for different density and viscosity ratios also increases relatively, as might be expected.

\begin{table}[t]
\centering
\footnotesize
\caption{Preconditioned GMRES iterations using two-phase PCD \eqref{GeneralisedPCD} / two-phase LSC \eqref{GeneralisedLSC} / LSC$_\text{D}$ \eqref{LSC_D} for the \mbox{$\ve{Q}_{2}$--$\,\ve{Q}_{1}$} solution to the steady lid driven cavity problem with density ratio $\widehat{\rho} = 1.2 \times 10^{-3}$, viscosity ratio $\widehat{\mu} = 1.8 \times 10^{-2}$ (air--water values), and varying Reynolds number $Re$ and grid size $h$.}
\label{TableGridSizeReynoldsNumber}
\tabulinesep=1.15mm
\begin{tabu}{|c|c|c|c|c|c|}
\hline
\multirow{2}{*}{$h$} & \multicolumn{5}{c|}{$Re$} \\
\cline{2-6}
& $10$ & $10^{1.5}$ & $100$ & $10^{2.5}$ & $1000$ \\
\hline
$\frac{1}{16}$ & 17 / 15 / 28 & 20 / 19 / 30 & 24 / 23 / 32 & 28 / 27 / 37 & 37 / 37 / 51 \\
\hline
$\frac{1}{32}$ & 19 / 18 / 54 & 21 / 20 / 57 & 25 / 28 / 63 & 29 / 34 / 73 & 35 / 30 / 89 \\
\hline
$\frac{1}{64}$ & 19 / 24 / 91 & 22 / 26 / 95 & \hphantom{0}27 / 41 / 103 & \hphantom{0}32 / 33 / 117 & \hphantom{0}36 / 46 / 148 \\
\hline
$\frac{1}{128}$ & \hphantom{0}20 / 48 / 145 & \hphantom{0}23 / 51 / 154 & \hphantom{0}28 / 56 / 166 & \hphantom{0}33 / 62 / 192 & \hphantom{0}38 / 67 / 244 \\
\hline
$\frac{1}{256}$ & \hphantom{0}20 / 69 / 243 & \hphantom{0}24 / 71 / 250 & \hphantom{0}29 / 81 / 267 & \hphantom{0}36 / 91 / 308 & \hphantom{0}44 / 99 / 396 \\
\hline
\end{tabu}
\end{table}

We mention here that our experience with the LSC variant known as LSC$_\text{D}$ \cite{MayAndMoresi,RehmanGeenenVuikSegalAndMacLachlan} shows a significant dependence on the grid size $h$. This is demonstrated in \Cref{TableGridSizeReynoldsNumber}, which allows a direct comparison between LSC$_\text{D}$ and our two-phase LSC approach. In particular, on the finest grid, LSC$_\text{D}$ requires more than triple the number of iterations required by our two-phase LSC approach. On finer meshes we expect this difference only to increase due to the poor scalability of LSC$_\text{D}$. We also note here that, while our two-phase PCD and, to a much lesser extent, LSC preconditioners show nearly mesh-independent behaviour, it may be that boundary adjustments must be made to achieve full mesh-independence, as is necessary for the original preconditioners \cite{ElmanSilvesterAndWathen}.

Before moving to unsteady flow, we provide analogous results for two-phase PCD when utilising our other choices of finite element pairings. \Cref{TableGridSizeReynoldsNumberPCDQ1-Q1AndQ2-P1} provides results for both (stabilised) \mbox{$\ve{Q}_{1}$--$\,\ve{Q}_{1}$} elements and \mbox{$\ve{Q}_{2}^{}$--$\,\ve{P}_{-1}^{}$} elements. 
We see similar behaviour with respect to mesh refinement and Reynolds number as with \mbox{$\ve{Q}_{2}$--$\,\ve{Q}_{1}$} (similar trends are also found in other studies not shown). The performance of two-phase PCD for \mbox{$\ve{Q}_{2}^{}$--$\,\ve{P}_{-1}^{}$} is particularly notable for being mesh-independent, aside from the large Reynolds number cases. We note that, in such large Reynolds number cases, the iteration counts initially decrease; this is attributed to the fact that on coarse meshes the flow is not well resolved or captured by the discrete operators used.

\begin{table}[t]
\centering
\footnotesize
\caption{Preconditioned GMRES iterations using two-phase PCD \eqref{GeneralisedPCD} for the \mbox{$\ve{Q}_{1}$--$\,\ve{Q}_{1}$} and \mbox{$\ve{Q}_{2}$--$\,\ve{P}_{-1}$} solution to the steady lid driven cavity problem with density ratio $\widehat{\rho} = 1.2 \times 10^{-3}$, viscosity ratio $\widehat{\mu} = 1.8 \times 10^{-2}$ (air--water values), and varying Reynolds number $Re$ and grid size $h$.}
\label{TableGridSizeReynoldsNumberPCDQ1-Q1AndQ2-P1}
\tabulinesep=1.15mm
\captionsetup{position=top}
\subfloat[$\ve{Q}_{1}$--$\,\ve{Q}_{1}$ elements]{
\label{TableGridSizeReynoldsNumberPCDQ1-Q1}
\begin{tabu}{|c|c|c|c|c|c|}
\hline
\multirow{2}{*}{$h$} & \multicolumn{5}{c|}{$Re$} \\
\cline{2-6}
& $10$ & $10^{1.5}$ & $100$ & $10^{2.5}$ & $1000$ \\
\hline
$\frac{1}{32}$  & 16 & 18 & 22 & 27 & 44 \\
\hline
$\frac{1}{64}$  & 16 & 19 & 24 & 29 & 34 \\
\hline
$\frac{1}{128}$ & 17 & 20 & 25 & 30 & 35 \\
\hline
$\frac{1}{256}$ & 17 & 20 & 25 & 32 & 37 \\
\hline
$\frac{1}{512}$ & \hphantom{0}17\hphantom{0} & \hphantom{0}20\hphantom{0} & \hphantom{0}26\hphantom{0} & \hphantom{0}33\hphantom{0} & \hphantom{0}43\hphantom{0} \\
\hline
\end{tabu} } \hspace*{\fill}
\subfloat[$\ve{Q}_{2}$--$\,\ve{P}_{-1}$ elements]{
\label{TableGridSizeReynoldsNumberPCDQ2-P1}
\begin{tabu}{|c|c|c|c|c|c|}
\hline
\multirow{2}{*}{$h$} & \multicolumn{5}{c|}{$Re$} \\
\cline{2-6}
& $10$ & $10^{1.5}$ & $100$ & $10^{2.5}$ & $1000$ \\
\hline
$\frac{1}{16}$  & 17 & 21 & 27 & 38 & 73 \\
\hline
$\frac{1}{32}$  & 18 & 21 & 27 & 34 & 55 \\
\hline
$\frac{1}{64}$  & 18 & 21 & 27 & 35 & 40 \\
\hline
$\frac{1}{128}$ & 18 & 21 & 27 & 34 & 42 \\
\hline
$\frac{1}{256}$ & \hphantom{0}18\hphantom{0} & \hphantom{0}21\hphantom{0} & \hphantom{0}27\hphantom{0} & \hphantom{0}36\hphantom{0} & \hphantom{0}49\hphantom{0} \\
\hline
\end{tabu} }
\end{table}

We now consider the time-dependent case ($\alpha = 1$). We remark that our numerical tests give similar dependence on the density and viscosity ratios, as well as on the grid size, to the steady case. The performance with Reynolds number is linked to the time-step, as illustrated in \Cref{TableTimeStepReynoldsNumber} using \mbox{$\ve{Q}_{2}$--$\,\ve{Q}_{1}$} elements. For a large time-step the behaviour is similar to that of the steady case. However, once the time-step is small enough, we see that the trend in Reynolds number flattens and then reverses so that higher Reynolds number flows require similar or less iterations for convergence. As anticipated, when solving with a smaller time-step fewer iterations are required since the time-stepping term becomes dominant, pushing the eigenvalues of the linear system away from zero and thus allowing faster convergence of the iterative solver.

\begin{table}[t]
\centering
\footnotesize
\caption{Preconditioned GMRES iterations using two-phase PCD \eqref{GeneralisedPCD} / two-phase LSC \eqref{GeneralisedLSC} for the \mbox{$\ve{Q}_{2}$--$\,\ve{Q}_{1}$} solution to the lid driven cavity problem with density ratio $\widehat{\rho} = 1.2 \times 10^{-3}$, viscosity ratio $\widehat{\mu} = 1.8 \times 10^{-2}$ (air--water values), $h = 1/128$, and varying Reynolds number $Re$ and time-step $\Delta t$.}
\label{TableTimeStepReynoldsNumber}
\tabulinesep=1.15mm
\begin{tabu}{|c|c|c|c|c|c|}
\hline
\multirow{2}{*}{$\Delta t$} & \multicolumn{5}{c|}{$Re$} \\
\cline{2-6}
& $10$ & $10^{1.5}$ & $100$ & $10^{2.5}$ & $1000$ \\
\hline
$10$ & 20 / 47 & 23 / 51 & 27 / 55 & 32 / 59 & 36 / 63 \\
\hline
$1$ & 19 / 46 & 21 / 47 & 23 / 45 & 24 / 45 & 25 / 41 \\
\hline
$\frac{1}{10}$ & \hphantom{0}16 / 39\hphantom{0} & \hphantom{0}16 / 37\hphantom{0} & \hphantom{0}16 / 35\hphantom{0} & \hphantom{0}15 / 32\hphantom{0} & \hphantom{0}16 / 28\hphantom{0} \\
\hline
\end{tabu}
\end{table}

The results of this section for a simplified two-phase test problem suggest that our new preconditioners perform reasonably well across a wide range of parameters and extend the utility of the PCD and LSC methodologies to two-phase flow. For the case of enclosed flow (without considering boundary adjustments) two-phase PCD appears favourable, showing nearly mesh-independent behaviour and only requiring a single Laplacian-type sub-problem to be solved per application of the Schur complement preconditioner. We now focus on the two-phase PCD approach. To gain further insight we consider a more challenging and realistic test problem within a framework which incorporates additional features pertinent to simulations of two-phase flow.

\section{Numerical results for a dynamic dam-break problem}
\label{sec:numericalresults2}

Free-surface models that accurately describe complicated air--water flow dynamics are an important application of the two-phase Navier--Stokes equations. For example, level set and volume-of-fluid methods can be combined with two-phase Navier--Stokes equations to simulate intricate hydraulic processes such as waves crashing into coastal structures. For these multi-physics models, the Navier--Stokes equations are just one component in a larger system of equations \cite{KeesEtAl}. In this section, we apply the two-phase PCD preconditioner in a free-surface model to gauge its effectiveness for dynamic multi-physics problems. In nearly all cases, solving the discrete linearised Navier--Stokes equations is the most time consuming part of a multi-physics simulation, so superior preconditioners can dramatically improve computational run times. As a result, the effectiveness of a preconditioner in a dynamic setting is very relevant to practitioners.

In the following study, we use the RANS2P module of the Proteus computational methods and simulation toolkit (\texttt{http://proteustoolkit.org}). RANS2P is a fluid dynamics software module developed at the U.S. Army Corps of Engineers for solving large air--water free-surface problems that arise in coastal and hydraulic applications. For a detailed description of the two-phase Navier--Stokes model within the RANS2P module, see \cite{ProteusTechnicalReport}.

As with other dynamic free-surface models, the two-phase Navier--Stokes model in RANS2P has several key characteristics that affect preconditioner performance. First, the dynamic nature of the model problems involved typically requires taking short time-steps. Consequently, the large temporal component of the problem diminishes the effect of the advective terms at each time-step, making the system easier to solve than a steady-state problem. A second feature of such free-surface models is that they often simulate high Reynolds number flows. As such, in order to produce meaningful results, we must modify the Navier--Stokes equations with numerical stabilisation terms. RANS2P uses a variation of ASGS stabilisation for the advection and pressure. Further details are available in \cite{ProteusTechnicalReport,Hughes,Codina,Tezduyar}.

Before moving forward, we say a few words about applying $\widehat{F} \, \! ^{-1}$ (recall \eqref{BlockPrecon}). For large problems, it is important that the action of $\widehat{F} \, \! ^{-1}$ is scalable and provides a reliable approximation to $F^{-1}$. Thus, developing effective methods to approximate $F \, \! ^{-1}$ is an active area of our research. For advection dominated flows, however, finding an effective and robust $\widehat{F} \, \! ^{-1}$ can be difficult \cite{OlshanskiiAndReuskenII}. Further complicating matters is the ASGS stabilisation used and the dynamic two-phase nature of the simulations. We have implemented an approach which treats $\widehat{F} \, \! ^{-1}$ as a block preconditioned GMRES iteration where the sub-blocks are approximated using multigrid methods and this has shown signs of being effective both in our work and others \cite{HeAndNeytcheva}; see also results in \cite{AxelssonHeAndNeytcheva}. Unfortunately, maintaining solver robustness throughout an entire multi-physics simulation remains an issue. As such, we defer discussion of this important topic to future work and focus our efforts here on the Schur complement approximation, taking $\widehat{F} \, \! ^{-1}$ to be given by the sparse direct solver SuperLU\_DIST \cite{SUPERLU,LiAndDemmel}.

We now comment on some implementation details of the PCD approximation used in RANS2P. Firstly, one component of the Navier--Stokes stabilisation is a shock capturing diffusion term. Thus, to better reflect the simulation viscosity, the viscosity used within the definition of two-phase PCD is the sum of the shock capturing term and the physical viscosity. Secondly, we use diagonal lumped mass matrices within two-phase PCD. Finally, even though the RANS2P module uses Newton iteration to solve the Navier--Stokes equations, rather than Picard iteration, the form of the PCD approximation does not change (see \cite{ElmanSilvesterAndWathen}). For additional details see \cite{ProteusTechnicalReport}.

As a basis for comparison, we consider the SIMPLE method (see \cite{ElmanSilvesterAndWathen,RehmanVuikAndSegal-SIMPLE}) where
\begin{align}
\label{SIMPLE}
\widehat{S}^{-1} = ( B \, \text{diag}(F)^{-1} B^{T} + C)^{-1}.
\end{align}
For steady problems, SIMPLE is not a competitive alternative to two-phase PCD since the $\text{diag}(F)$ term does not capture well the advective features of the flow. However, for dynamic problems with short time-steps, the method is effective in many cases.

Neither SIMPLE nor the two-phase PCD method offers a clear implementation advantage over the other. Indeed, at each time-step, three (sparse) discrete operators must be constructed for two-phase PCD. In contrast, SIMPLE \eqref{SIMPLE} is algebraic and constructed using blocks taken directly from the global linear system. For sufficiently large problems, however, the matrix--matrix product used in the construction of \eqref{SIMPLE} is expensive. In terms of applying the approximations, two-phase PCD requires one V-cycle to solve the pressure Laplacian-type term $A_{p}^{(1/\rho)}$, a vector--scalar product, two diagonal matrix--vector products, a sparse matrix--vector product, and two vector sums. By comparison, applying SIMPLE only requires a single multigrid V-cycle. Our experience, however, suggests that SIMPLE requires a more accurate multigrid solver than two-phase PCD to serve as a reliable preconditioner in the RANS2P module. In each approach we make use of BoomerAMG from HYPRE \cite{FalgoutandYang} for the V-cycle.

In the following examples we use right preconditioned GMRES to compare the SIMPLE and two-phase PCD methods. This ensures the convergence of the linear solver is measured using non-preconditioned residuals. In RANS2P, our experience suggests that non-preconditioned residuals produce more stable results and provide a reliable basis for comparing different preconditioners. 

For dynamic free-surface simulations, the linear solver termination criterion used depends on the nonlinear solver tolerance. We always use a relative residual tolerance of $10^{-6}$ for the linear solver. However, we also include an absolute linear solver tolerance that is one order of magnitude smaller than the absolute nonlinear solver tolerance. For example, if the absolute tolerance of the nonlinear solver is $10^{-7}$, then the linear solver will include an absolute tolerance of $10^{-8}$ for termination. This approach maintains nonlinear solver convergence while limiting unnecessary work in the linear solver.

Before examining preconditioner performance for a complete multi-physics model, we first present results using (unstructured) \mbox{$\ve{P}_{1}$--$\,\ve{P}_{1}$} elements for the steady lid driven cavity problem using the RANS2P module. Note that we measure the mesh refinement here using $h = \sqrt{2 T_{A}}$, where $T_{A}$ is the area of the largest triangle. \Cref{TableGridSizeReynoldsNumberProteus} shows analogous results to \Cref{TableGridSizeReynoldsNumber} for the two-phase PCD and SIMPLE methods. We see that the SIMPLE method fails to produce scalable results, with iteration counts increasing as we refine the mesh. In contrast, the two-phase PCD approach exhibits mesh-independent behaviour similar to the \mbox{$\ve{Q}_{1}$--$\,\ve{Q}_{1}$} element results from IFISS in \Cref{TableGridSizeReynoldsNumberPCDQ1-Q1}. It is also worth mentioning that the iteration level of the \mbox{$\ve{P}_{1}$--$\,\ve{P}_{1}$} elements is higher than the \mbox{$\ve{Q}_{1}$--$\,\ve{Q}_{1}$} elements because different stabilisation methods are used.

\begin{table}[t]
\centering
\footnotesize
\caption{Preconditioned GMRES iterations using two-phase PCD \eqref{GeneralisedPCD} / SIMPLE \eqref{SIMPLE} for the \mbox{$\ve{P}_{1}$--$\,\ve{P}_{1}$} solution to the steady lid driven cavity problem with density ratio $\widehat{\rho} = 1.2 \times 10^{-3}$, viscosity ratio $\widehat{\mu} = 1.8 \times 10^{-2}$ (air--water values), and varying Reynolds number $Re$ and mesh size $h$. The total number of degrees of freedom (DOFs) for each problem size is also listed for comparison with \Cref{TableNDOFs}.}
\label{TableGridSizeReynoldsNumberProteus}
\tabulinesep=1.15mm
\begin{tabu}{|c|c|c|c|c|c|c|}
\hline
\multirow{2}{*}{$h$} & \multirow{2}{*}{DOFs} & \multicolumn{5}{c|}{$Re$} \\
\cline{3-7}
& & $10$ & $10^{1.5}$ & $100$ & $10^{2.5}$ & $1000$ \\
\hline
$\frac{1}{16}$ & $5,\!115$ & \hphantom{0}25 / 18\hphantom{0} & \hphantom{0}27 / 18\hphantom{0} & \hphantom{0}29 / 23\hphantom{0} & \hphantom{0}36 / 26\hphantom{0} & \hphantom{0}49 / 25\hphantom{0} \\
\hline
$\frac{1}{32}$ & $20,\!085$ & \hphantom{0}24 / 38\hphantom{0} & \hphantom{0}26 / 23\hphantom{0} & \hphantom{0}29 / 32\hphantom{0} & \hphantom{0}35 / 36\hphantom{0} & \hphantom{0}47 / 39\hphantom{0} \\
\hline
$\frac{1}{64}$ & $79,\!152$ & 25 / 56 & 27 / 55 & 29 / 56 & 35 / 45 & 46 / 56 \\
\hline
$\frac{1}{128}$ & $308,\!613$ & 24 / 53 & 27 / 85 & 29 / 77 & 34 / 83 & 46 / 84 \\
\hline
$\frac{1}{256}$ & $1,\!191,\!810$ & \hphantom{0}24 / 116 & \hphantom{0}27 / 116 & \hphantom{0}30 / 105 & \hphantom{0}35 / 112 & \hphantom{0}47 / 115 \\
\hline
\end{tabu}
\end{table}

\begin{figure}[t]
\centering
\subfloat[$t=0$]{ \includegraphics[width=.488\linewidth]{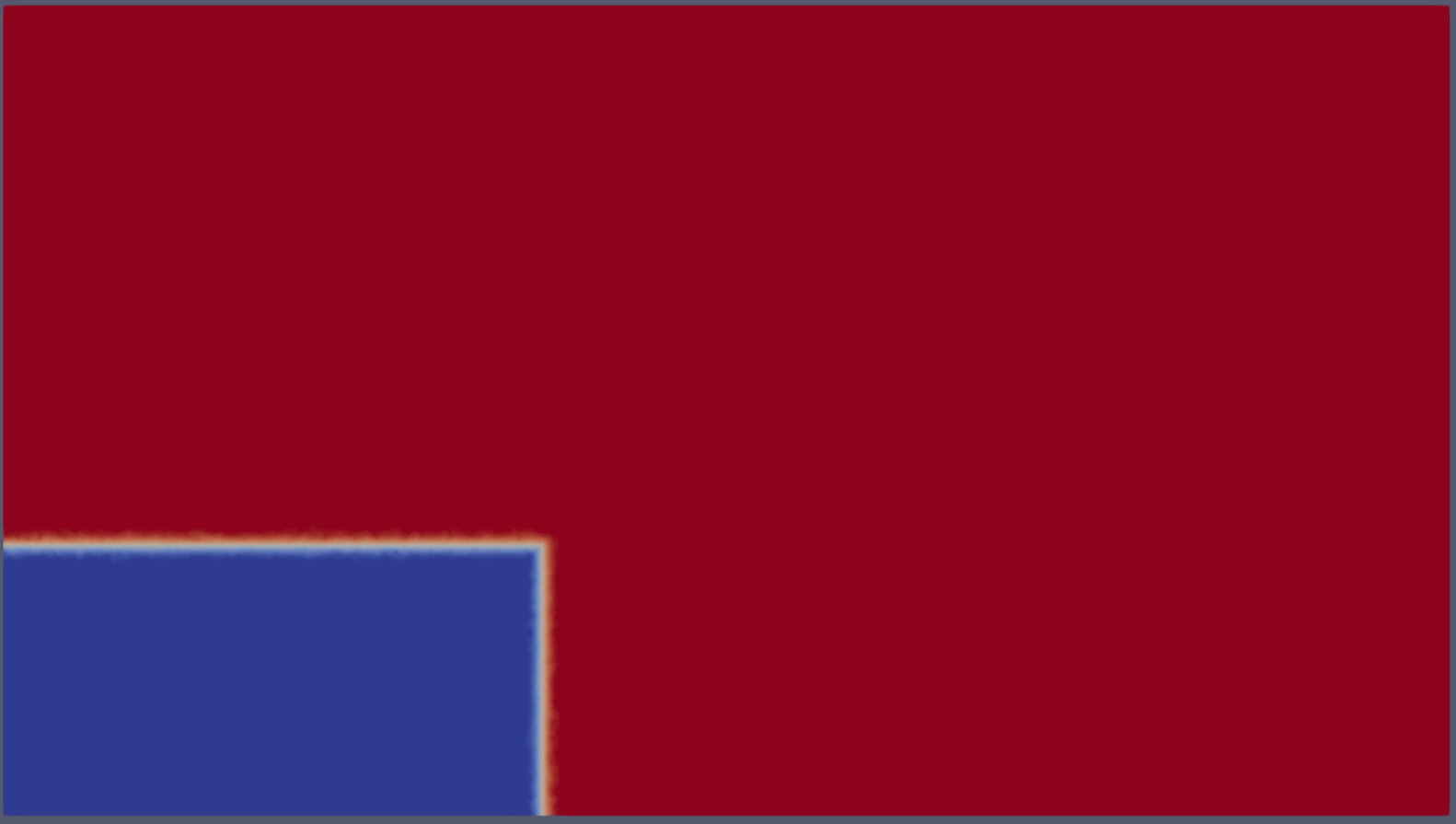}} \hspace*{\fill}
\subfloat[$t=0.5$]{ \includegraphics[width=.488\linewidth]{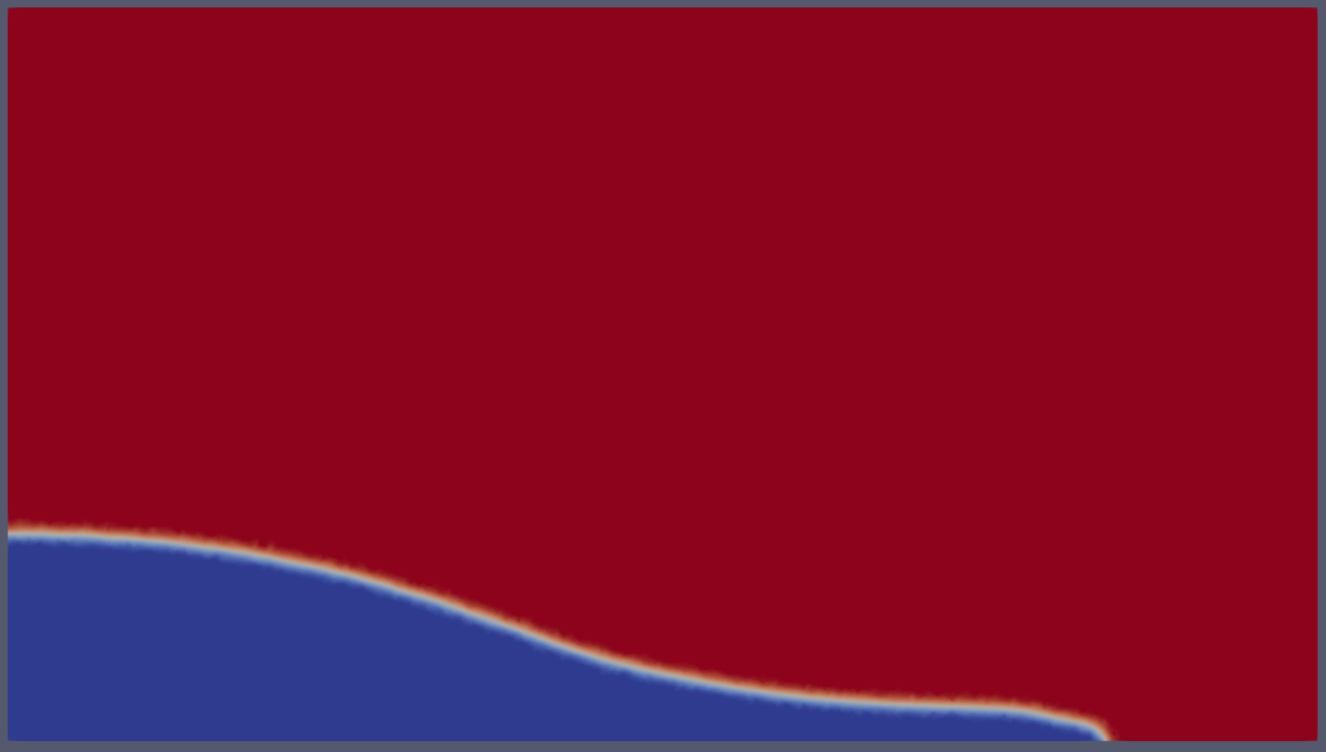}}
\\
\subfloat[$t=1$]{ \includegraphics[width=.488\linewidth]{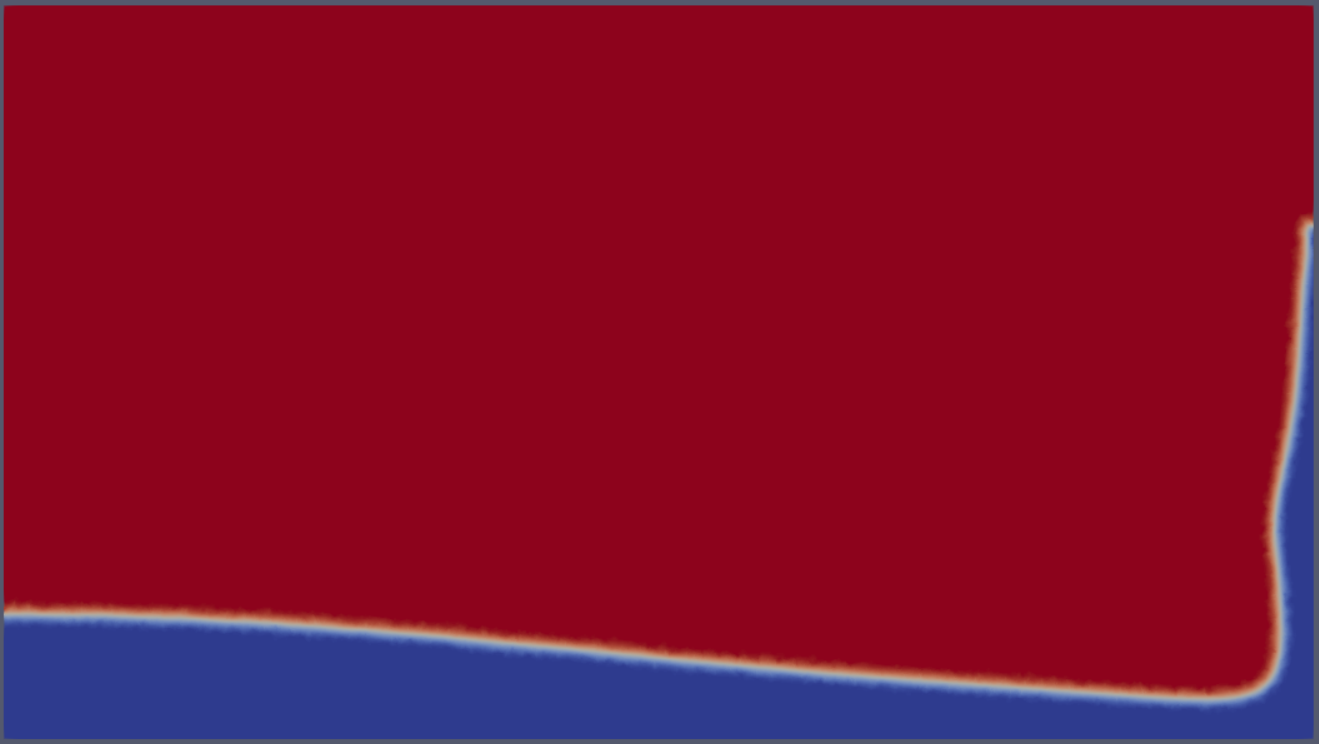}} \hspace*{\fill}
\subfloat[$t=1.5$]{ \includegraphics[width=.488\linewidth]{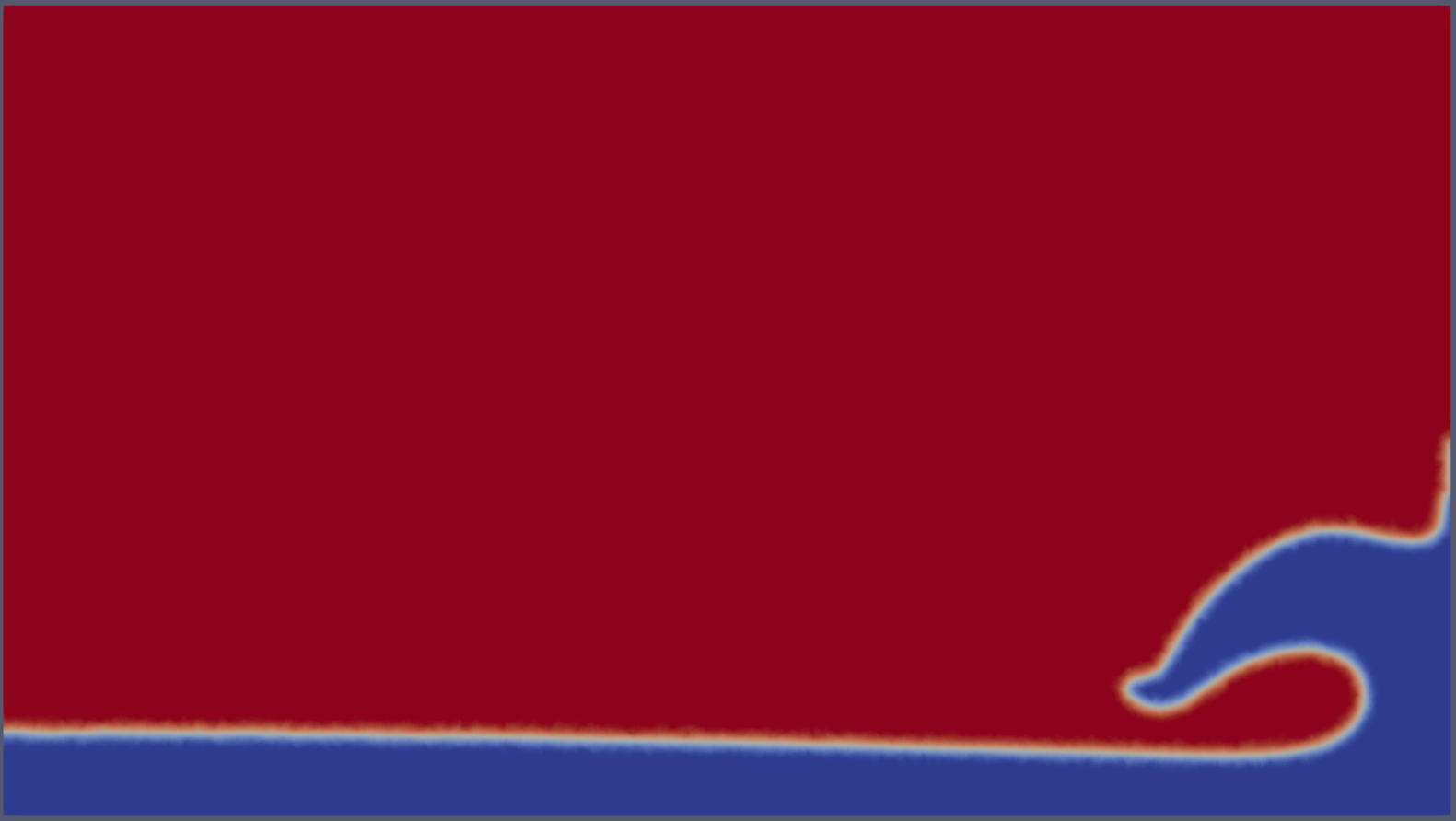}}
\\
\subfloat[$t=1.75$]{ \includegraphics[width=.488\linewidth]{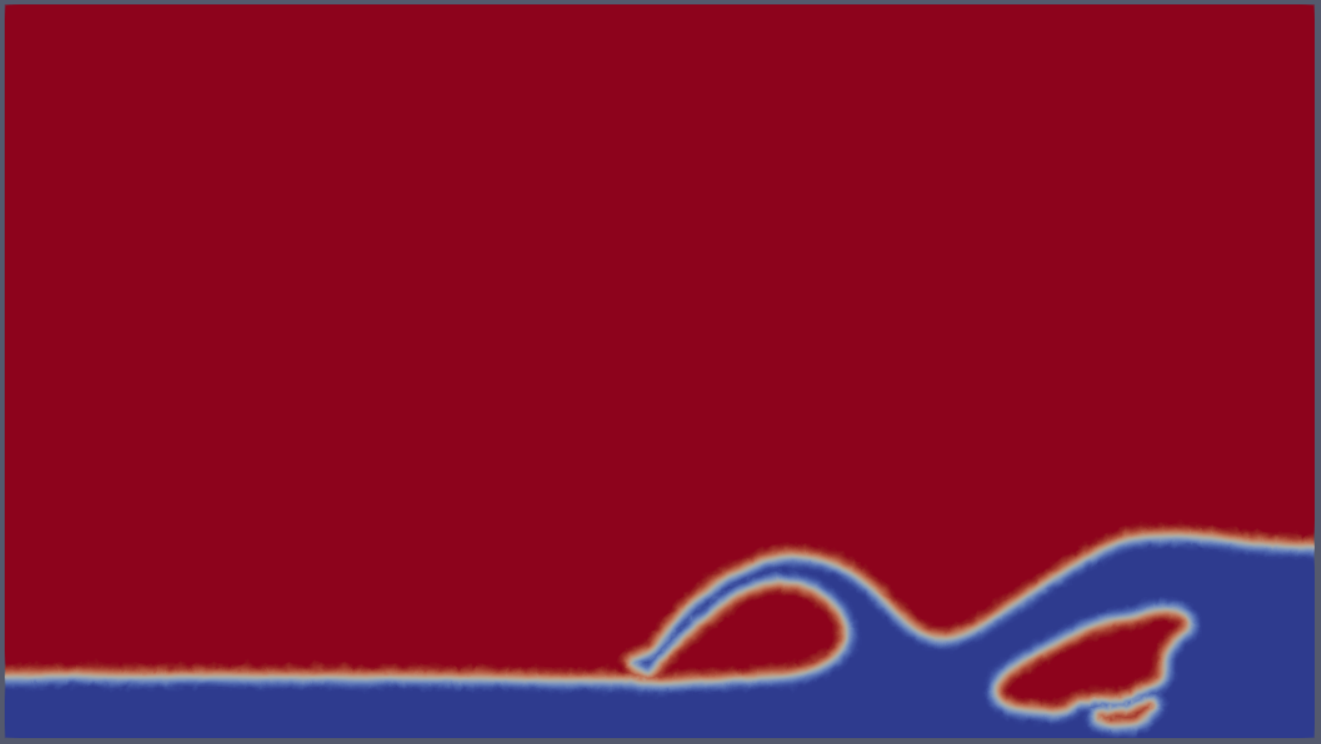}} \hspace*{\fill}
\subfloat[$t=2$]{ \includegraphics[width=.488\linewidth]{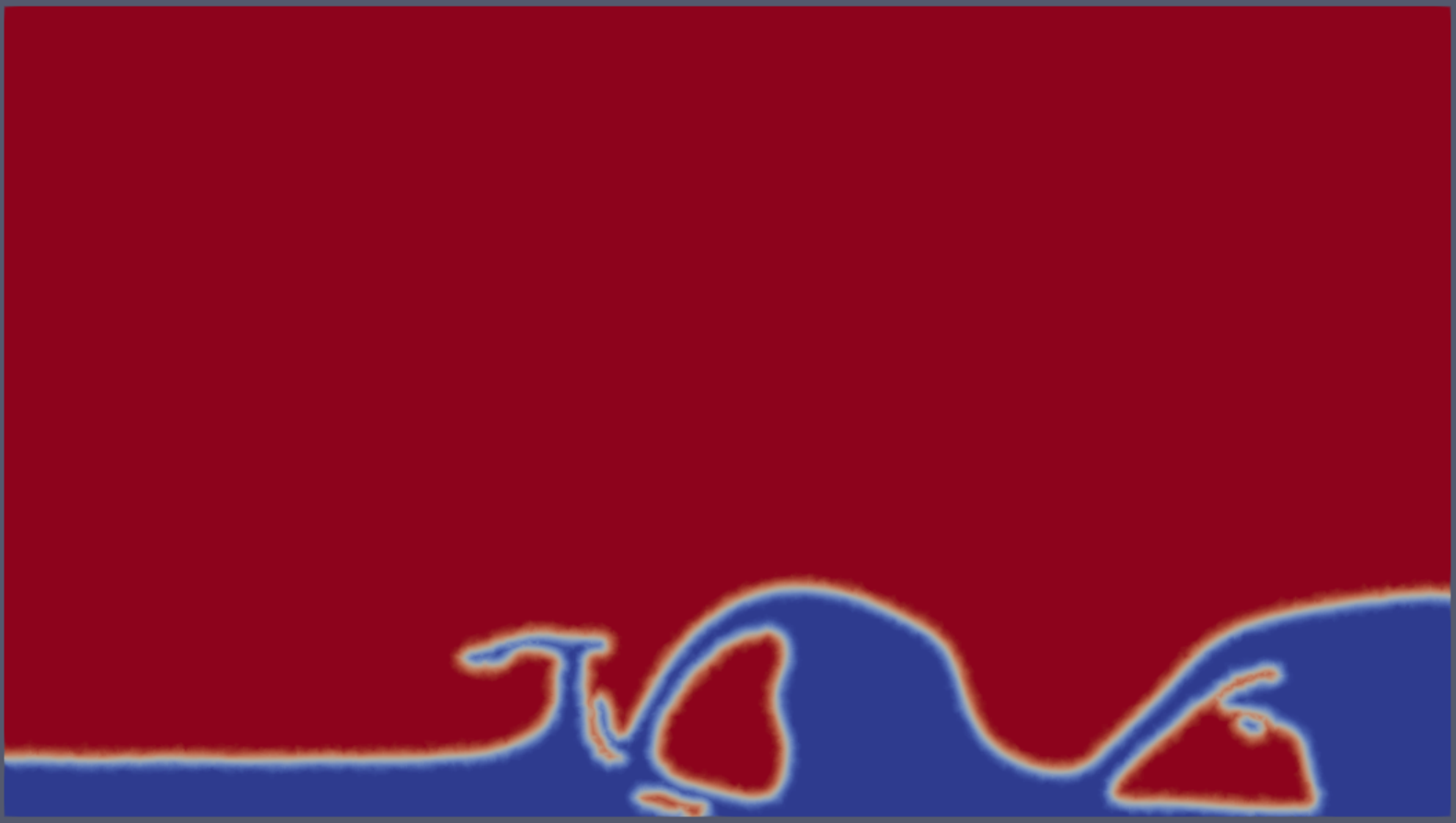}}
\caption{Evolution of the dam-break simulation in Proteus at selected points in time. The VOF (volume-of-fluid) is plotted with blue representing the water phase and red being the air.}
\label{FigureDam-breakSimulation}
\end{figure}

Next we consider the 2D dam-break benchmark problem described in \cite{CollagrossiAndLandrini,ZhouDeKatAndBuchner}, which requires using the complete multi-physics RANS2P model. The domain is rectangular, with $\Omega = (0, 3.22) \times (0, 1.8)$, and the free-slip condition \eqref{BoundaryConditionsN} is applied everywhere on the boundary $\partial \Omega$. Initially, there is a standing column of water in $\Omega_{1} = (0, 1.2) \times (0, 0.6)$ with the remaining space being air. The simulation runs for a two second time interval and begins as the column of water collapses under gravity and proceeds to collide with the right-hand wall of the tank. This collision creates a wave and ultimately topological changes in the phases. \Cref{FigureDam-breakSimulation} displays several snapshots of the simulation. The dam-break problem provides a good benchmark for testing the two-phase PCD and SIMPLE preconditioners because its features are typical of many dynamic, multi-physics problems of practical interest.

To compare the scaling performance of two-phase PCD and SIMPLE, we consider two simulations. In the first, time-steps are selected to ensure that the CFL number is less than or equal to 0.9. Such restrictions are often necessary for nonlinear solver convergence and solution accuracy. However, in some important cases this restriction on the CFL number is not strictly necessary.\footnote{Accurate computation of relevant quantities of interest, such as drag force, for fixed hydraulic structures or vessels, frequently results in quasi-steady flows. In particular, the free surface may tend towards a steady wake structure or standing wave pattern, and this structure dominates the force on the given structure. In these cases, it is frequently desirable to use a fixed time-step that results in CFL numbers significantly larger than one. Time-stepping is then carried out until the quasi-steady hydrodynamic conditions are reached or the quantity of interest has reached a constant or steady periodic value.} Thus, in the second simulation, we use a fixed time-step of $\Delta t = 0.01$. In this case, the CFL number is larger than one for much of the simulation, reaching a maximum of $20.5$ and typically being above $2.5$. Nonetheless, the time-step is still small enough to achieve nonlinear solver convergence and solution accuracy. For both simulations, we analyse the average and maximum number of GMRES iterations required at five different levels of mesh refinement. These mesh refinement levels are selected so that, by the final refinement, the physics of the simulation is sufficiently resolved to perform relevant engineering analysis. The dam-break timings were collected using 8 cores of a dedicated 2.3-GHz Intel Xeon Haswell processor with 128 GBytes of DDR4 memory on the Topaz supercomputer in the Department of Defense High Performance Computing Modernization Program.

\begin{figure}[t]
\centering
\subfloat[Two-phase PCD]{ \includegraphics[width=.488\linewidth,clip,trim=2mm -1cm 3cm 0]{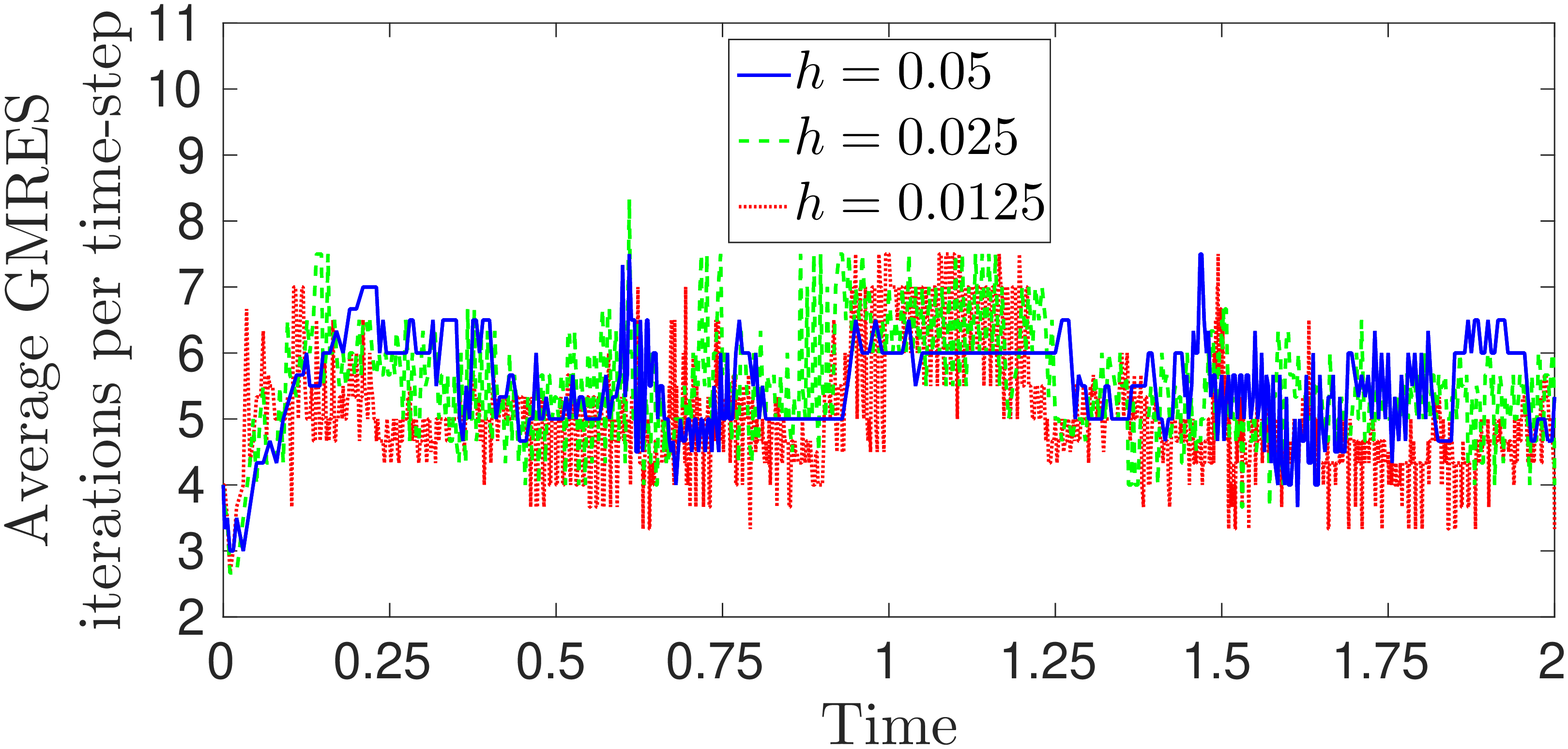}} \hspace*{\fill}
\subfloat[SIMPLE]{ \includegraphics[width=.488\linewidth,clip,trim=2mm -1cm 3cm 0]{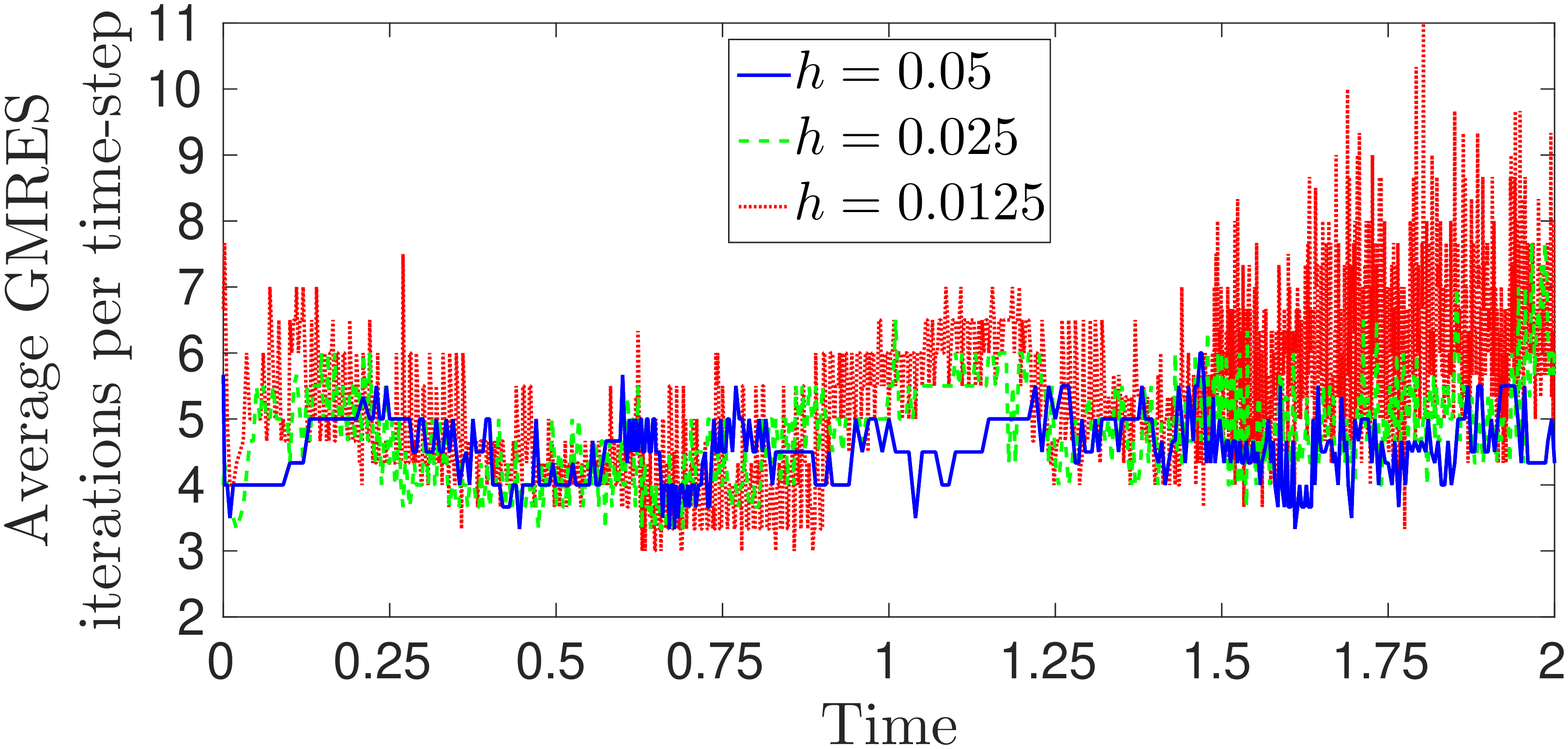}}
\caption{Average preconditioned GMRES iterations per time-step for the two-phase PCD and SIMPLE preconditioners when the CFL number is less than or equal to $0.9$.}
\label{fig:pcd_plot_cfl_leq_09}
\end{figure}

\Cref{TableFixedCFL} presents the average and maximum number of GMRES iterations taken during the first simulation with a restricted CFL number. These results suggest that, on average, the SIMPLE and two-phase PCD preconditioners both scale well with the mesh size. However, \Cref{TableFixedCFL} also reveals that the maximum number of iterations required by the SIMPLE method increases as the mesh is refined. As seen in \Cref{fig:pcd_plot_cfl_leq_09}, the increase in maximum iterations of the SIMPLE preconditioner occurs as the air and water phases begin to undergo topological changes around one and a half seconds into the simulation. Indeed, as the water phase reconnects with itself, it generates a pressure that causes the air phase to accelerate and increase the advective features of the simulation. Since SIMPLE only uses the diagonal elements of the matrix $F$, it appears unable to fully capture these additional advection dynamics. In contrast, the two-phase PCD preconditioner scales well during this mixing phase of the simulation.

\begin{table}[t]
\centering
\footnotesize
\caption{The average / maximum number of GMRES iterations and simulation run times (in minutes) required across different meshes when running the dam-break problem ensuring the CFL number is less than or equal to $0.9$.}
\label{TableFixedCFL}
\tabulinesep=1.15mm
\begin{tabu}{|c|c|c|c|c|c|}
\hline & $h=0.2$ & $h=0.1$ & $h=0.05$ & $h=0.025$ & $h=0.0125$ \\
\hline
Two-phase PCD & 5 / 8\hphantom{0} (0.5) & 5 / 9\hphantom{0} (0.6) & 5 / 10 (2.6) & 5 / 11 (14.7) & 5 / 10 (126.0) \\
\hline
SIMPLE & 4 / 10 (0.4) & 4 / 10 (0.6) & 4 / 10 (2.4) & 5 / 13 (15.5) & 5 / 17 (140.1) \\
\hline
\end{tabu}
\end{table}

\Cref{TableFixedCFL} also reveals that, on the most refined meshes, the two-phase PCD method is faster than the SIMPLE approach. One reason for this is that the two-phase PCD preconditioner requires fewer GMRES iterations than the SIMPLE preconditioner during the mixing phase of the simulation. A second reason is that the AMG method applied within the SIMPLE preconditioner requires more computational effort than the AMG method used in the two-phase PCD preconditioner. Finally, the two-phase PCD method tends to exceed the linear solver threshold by a larger margin than the SIMPLE approach, leading to slightly smaller residual norms in the nonlinear solver. Interestingly, this difference slightly reduces the computational effort needed to solve other components in the full RANS2P model.

Results for the second simulation, using a fixed time-step $\Delta t = 0.01$, are shown in \Cref{TableFixedTimeStep1} and suggest that, on coarse meshes, SIMPLE and the two-phase PCD preconditioner are competitive with one another. In contrast to the first simulation, however, as the mesh is refined, the performance of the SIMPLE method rapidly deteriorates while the two-phase PCD preconditioner remains relatively stable, with iteration counts increasing only modestly. Again, two-phase PCD is preferable based on the timing results. These results are consistent with the steady-state performance observed above in \Cref{TableGridSizeReynoldsNumberProteus}. As the mesh is refined for a fixed time-step, the advective features of the system become more pronounced and the CFL number increases. As observed for the steady lid driven cavity problem, the SIMPLE approach does not capture the features of an advection dominated flow well enough to provide a robust preconditioner. The two-phase PCD preconditioner, however, does account for such features and thus remains capable of producing stable, reliable results in this setting.

\begin{table}[t]
\centering
\footnotesize
\caption{The average / maximum number of GMRES iterations and simulation run times (in minutes) required across different meshes when running the dam-break problem with fixed $\Delta t = 0.01$. Note that $\times$ indicates the simulation stopped due to failure in the solver convergence.}
\label{TableFixedTimeStep1}
\tabulinesep=1.15mm
\begin{tabu}{|c|c|c|c|c|c|}
\hline 
 & $h=0.2$ & $h=0.1$ & $h=0.05$ & $h=0.025$ & $h=0.0125$ \\
\hline
Two-phase PCD & 4 / 8\hphantom{0} (0.5) & 5 / 9\hphantom{0} (0.6) & 8 / 14 (1.6) & 11 / 25 (5.8) & 14 / 34 (26.5) \\
\hline
SIMPLE & 4 / 10 (0.4) & 4 / 10 (0.6) & 5 / 10 (1.5) & 10 / 32 (6.6) & $\times$ \\
\hline
\end{tabu}
\end{table}

Overall, our results for the two-phase PCD preconditioner in a free-surface, multi-physics setting are encouraging. When a restricted CFL number is used, the two-phase PCD preconditioner slightly outperforms the SIMPLE method both in terms of the reducing the number of GMRES iterations required as well as delivering faster run times. As the CFL number of the flow increases, two-phase PCD demonstrates a significant improvement over the SIMPLE method due to its superior steady-state performance. Together, these results suggest that the two-phase PCD approach can be effectively used as an approximation to the inverse Schur complement in coupled free-surface problems.

\section{Conclusions}
\label{sec:conclusions}

The application of PCD and LSC preconditioning techniques for the Navier--Stokes equations has proved to be an effective approach for computing flows of a single Newtonian phase. In this work we have generalised the formulation of these preconditioners to the case of the variable coefficient Navier--Stokes equations arising in models of two-phase incompressible flow. For the PCD method this requires a new form of the preconditioner, a point made clear through the relation of PCD with the Cahouet--Chabard preconditioner. While a variant of LSC has previously been proposed for variable viscosity Stokes flow, we present a superior method for Navier--Stokes flow which directly appeals to the relevant commutator. Our new two-phase PCD and LSC preconditioners retain the favourable features of the original preconditioners and exhibit similar performance in our tests on a simplified enclosed flow problem. As for the original methods, additional consideration of boundary conditions will be necessary to gain fully mesh-independent behaviour, in particular for the case of flow problems with inflow and outflow boundary conditions when using the two-phase PCD preconditioner \cite{ElmanSilvesterAndWathen}. We intend to investigate this in future work.

Our results for a fully dynamic dam-break test problem showed that the two-phase PCD approach can be effectively used in coupled free-surface problems. Moreover, since it captures sufficient information about the advective term, it provides a stable preconditioner when larger time-steps are permissible. Nonetheless, further work is needed to investigate whether the additional terms added to the discrete equations can be incorporated within the two-phase PCD preconditioner to improve performance. Moreover, an important part of the block preconditioners which we consider here is a scalable approximation to the velocity block $F$. This topic requires further attention for variable coefficient Navier--Stokes equations and is an active area of investigation. Additionally, in practical applications, simulations of two-phase flow run in a high performance computing environment with a parallel implementation; it remains to be seen how efficient the preconditioners we have proposed are within this framework. We are currently exploring the parallel scalability of our approach on larger problems of interest. 

Finally, we note that we did not require the two-phase nature explicitly in the construction of the preconditioners, thus the techniques proposed here might also be investigated for their utility when more general forms of variable density and variable viscosity flows are proposed.

\section{Acknowledgements}
\label{sec:acknowledgements}

The authors would like to thank the associate editor and two anonymous referees for their careful reading and constructive comments which helped improve the content and presentation of this work.

\bibliographystyle{siamplain}
\bibliography{references}

\end{document}